\documentclass[10pt]{amsart}

\usepackage{tikz}
\usepackage{amssymb,amsmath,latexsym,graphicx}
\usepackage{charter,eucal}
\usepackage{amssymb}
\usepackage{amscd}
\usepackage[all,cmtip]{xy}
\usepackage{rotating}
\usetikzlibrary{decorations.markings}

\newtheorem{theorem}{Theorem }[section]

\newtheorem{lemma}[theorem]{Lemma}
\newtheorem{observation}[theorem]{Observation}

\newtheorem{remark}[theorem]{Remark}
\newtheorem{corollary}[theorem]{Corollary}
\newtheorem{proposition}[theorem]{Proposition}
\newtheorem{principle}[theorem]{\textsc{Principle}}

\newcommand{\bt}{\begin{theorem}}
\newcommand{\et}{\end{theorem}}
\newcommand{\bmt}{\begin{maintheorem}}
\newcommand{\emt}{\end{maintheorem}}
\newcommand{\bc}{\begin{corollary}}
\newcommand{\bl}{\begin{lemma}}
\newcommand{\ec}{\end{corollary}}
\newcommand{\el}{\end{lemma}}
\newcommand{\bo}{\begin{observation}}
\newcommand{\eo}{\end{observation}}
\newcommand{\bp}{\begin{proposition}}
\newcommand{\ep}{\end{proposition}}
\newcommand{\br}{\begin{remark}}
\newcommand{\er}{\end{remark}}
\newcommand{\bpr}{\begin{principle}}
\newcommand{\epr}{\end{principle}}

\def\Sym{\mathrm{Sym}}
\def\hol{\mathrm{hol}}
\def\soc{\mathrm{soc}}
\def\Out{\mathrm{Out}}
\def\Aut{\mathrm{Aut}}
\def\I{\mathop{\mathrm{I}}}

\def\PG{\mathbf{PG}}
\def\PSL{\mathbf{PSL}}
\def\GL{\mathbf{GL}}
\def\C{\mathbb{C}}
\def\AG{\mathbf{AG}}
\def\AGL{\mathbf{AGL}}

\def\eop{\hspace*{\fill}$\blacksquare$}
\def\eopr{\hspace*{\fill}$\nabla$}

\def\O{\mathbf{O}}
\def\Inn{\mathbf{Inn}}
\def\id{\mathrm{id}}
\def\bA{\mathcal{A}}

\newcommand{\PGL}{\mathbf{PGL}}
\newcommand{\F}{\mathbb{F}}

\newcommand{\mX}{\mathcal{X}}
\newcommand{\T}{\mathbf{T}}
\newcommand{\mA}{\mathcal{A}}

\newcommand{\mE}{\mathcal{E}}

\newcommand{\M}{\mathcal{M}}
\newcommand{\K}{\mathbb{K}}

\newcommand{\mS}{\mathcal{S}}
\newcommand{\mD}{\mathcal{D}}

\title[Isospectral drums and simple groups]{Isospectral drums and simple groups}


\author{Koen Thas}

\thanks{}

\address{Ghent University, Department of Mathematics, Krijgslaan 281, S25, B-9000 Ghent, Belgium}

\email{koen.thas@gmail.com}

\date{}

\begin{document}

\maketitle

\begin{abstract}
Virtually every known pair of isospectral but nonisometric manifolds | with as most famous members isospectral  bounded $\mathbb{R}$-planar domains which makes one ``not hear the shape of a drum'' \cite{Kac} | arise from the (group theoretical) Gassman-Sunada method. Moreover, all the known $\mathbb{R}$-planar examples (so counter examples to Kac's question) are constructed through a famous specialization of this method, called {\em transplantation}.\\
We first describe a number of very general classes of {\em length equivalent} manifolds, with as particular cases {\em isospectral} manifolds, in each of the constructions starting from a  given example that arises itself from the Gassman-Sunada method. 
The constructions include the examples arising from the transplantation technique (and thus in particular the planar examples). 
To that end, we introduce four properties | called FF, MAX, PAIR and INV | inspired by natural physical properties (which rule out trivial constructions), that are satisfied for each of the known  planar examples. \\
Vice versa, we show that length equivalent manifolds with FF, MAX, PAIR and INV which arise from the Gassman-Sunada method,  {\em must} fall under one of our prior constructions, thus describing a precise classification of these objects. \\
Due to the nature of our constructions and properties, a deep connection with finite simple groups occurs which seems, perhaps, rather surprising in the context of this paper. On the other hand, our properties define physically irreducible pairs of length equivalent manifolds | ``atoms'' of {\em general}
pairs of length equivalent manifolds, in that such a general pair of manifolds is patched up out of irreducible pairs | and that is precisely what simple groups are for general groups.  
%Contrary to the fact that most known examples arise from two-step nilpotent $(p)$-groups, our examples live at the complete other end of the %spectrum: they arise from finite simple groups.
%(SOLVABLE, 2-STEP NILP., NOT HERE)
\end{abstract}

\setcounter{tocdepth}{1}
\bigskip
%{\footnotesize
\tableofcontents
%}

\newpage
\section{Introduction and statement of the main results}

In 1966, in a celebrated paper \cite{Kac}, Mark Kac
formulated the famous question ``Can one hear the shape of a
drum?'': is it possible to find two  nonisometric Euclidean simply connected domains for which
the sets $\{ E_n \vert n \in\mathbb{N}\}$ of solutions of 
\begin{equation}
\label{helmholtz}
\Delta f+E f=0,
\end{equation}
where $\Delta$ is the $d$-dimensional Laplacian
and 
with $\Psi_{\vert \mbox{Boundary}} = 0$, are
identical? It was known very early, from
Weyl's formula, that one can ``hear'' the area of a drum and the
length of its perimeter. But could the
shape itself be retrieved from the spectrum?\\

Formally, answering ``no'' to Kac's question amounts to finding
{\em isospectral billiards}, that is, nonisometric billiards having
exactly the same eigenvalue spectrum (with multiplicities). 
Since the appearance of \cite{Kac},  an incredible amount of variations on
the theme have found their way through literature. 

Very early examples of flat tori sharing the same
eigenvalue spectrum had been found in 1964 by Milnor in $\mathbb{R}^{16}$ in his classic note \cite{Milnor}, from
nonisometric lattices of rank $16$ in $\mathbb{R}^{16}$.  Other examples of isospectral Riemannian manifolds
were constructed later, e.g., on lens spaces \cite{Ike} or on surfaces
with constant negative curvature \cite{Vig}. In 1982,  Urakawa produced the
first examples of isospectral domains in $\mathbb{R}^{n}$, $n\geq 4$ \cite{Ura}, applying work
of B\'{e}rard and Besson \cite{BerBes}.

In the late eighties, various other papers appeared, giving necessary conditions that any family of
billiards sharing the same spectrum should fulfill (\cite{Mel},
\cite{OsgPhiSar}, \cite{OsgPhiSar2}), and necessary conditions given as
inequalities on the eigenvalues were reviewed in \cite{Pro}.\\

\subsection{Planar counter examples}

\begin{figure}[ht]
\begin{center}
\includegraphics[width=10cm]{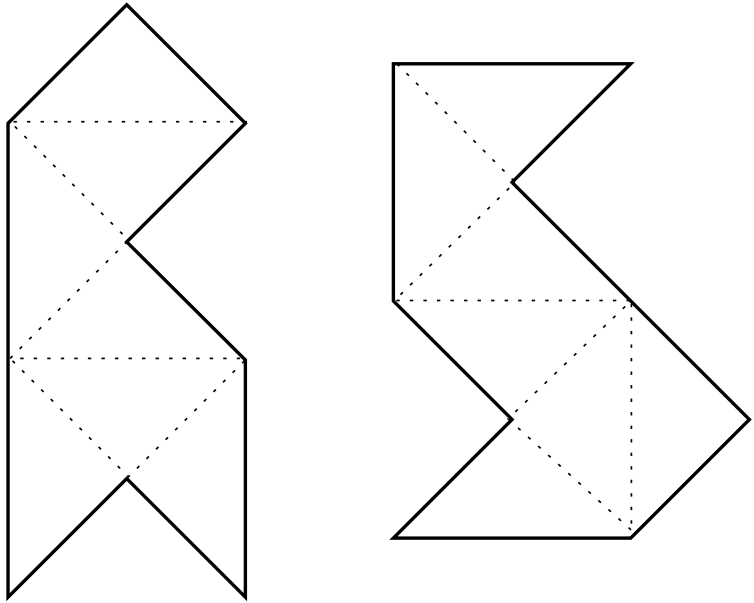}
\end{center}
\caption{The famous pair of isospectral billiards with seven
  half-square shaped base tiles. The dotted lines are just an eyeguide.
\label{celebrated}}
\end{figure}

It took almost 30 years after Kac's paper that the first example of
two-dimensional billiards having exactly the same spectrum was 
exhibited (in 1992). The pair was found by Gordon, Webb and  Wolpert in the paper
``Isospectral plane domains and surfaces via Riemannian orbifolds''
\cite{GWW}. The authors thus gave a ``no'' as a final answer to Kac's
question, and as a reply to Kac's paper, they published a paper titled ``One cannot hear the shape of a
drum'' \cite{GWW2}. The most popularized example is shown in Figure \ref{celebrated}.

Crucial for finding the example was a theorem by Sunada
 asserting that when two subgroups that are
``almost conjugate'' in a group that acts by isometries on a Riemannian manifold, the quotient
manifolds are isospectral. In fact, the other examples which were 
constructed after 1992 virtually all used Sunada's method. 
Later, the so-called ``transplantation technique'' was introduced, giving
an easier way of detecting isospectrality of planar billiards. (Still,
essentially only seventeen families of examples that say ``no'' to Kac's
question were constructed in a forty year period.) 

We refer to \cite{GirTha} for an extensive review on the modern theory after Kac.\\

\subsection{More on Sunada's method}

For $\M$ a Riemannian manifold, one defines the zeta function
\begin{equation}  \zeta_{\M}(s) = \sum_{i = 1}^{\infty}\lambda_i^{-s}, \ \ \texttt{Re}(s) \geq 0,                       \end{equation}
where $0 < \lambda_1 \leq \lambda_2 \leq \cdots$
are the nonzero eigenvalues of the Laplacian for $\M$.
The function $\zeta_{\M}$ has an analytic continuation to the whole plane, and it is well known that $\zeta_{\M_1}(s) = \zeta_{\M_2}(s)$ if and only if $\M_1$ and $\M_2$ are isospectral.

The next theorem is important.

\begin{theorem}[T.~Sunada \cite{Su}]
\label{sunada}
Let $\pi: \M \mapsto \M_0$ be a normal finite Riemannian covering with covering transformation group $G$, and let $\pi_1: \M_1 \mapsto \M_0$ and $\pi_2: \M_2 \mapsto \M_0$ be the coverings corresponding to the subgroups $H_1$ and $H_2$ of $G$, respectively. If the triplet $(G,H_1,H_2)$ satisfies the property that each conjugacy class of $G$ meets $H_1$ and $H_2$ in the same number of elements, then the zeta functions $\zeta_{\M_1}(s)$ and $\zeta_{\M_2}(s)$ are identical.
\end{theorem}

Triples as in the theorem are called ``Gassman-Sunada triples,'' and one also says that $H_1$ and $H_2$ are ``almost conjugate' in $G$ (and 
$(G,H_1,H_2)$ is an almost conjugate or AC-triple). 

Expressed in terms of isospectrality, we have the following version of Theorem \ref{sunada}.
Let $\M$ be a closed hyperbolic $n$-manifold, $G$ a group, and $H$ and $K$ almost conjugate subgroups of $G$. Then if $\pi_1(\M)$ admits a 
homomorphism onto $G$, the covers $\M_H$ and $\M_K$ associated to the pullback subgroups of $H$ and $K$ are isospectral.

This bridge between spectral theory and Group Theory gave a strong impetus to the construction theory for isospectral manifolds, and indeed,
the first counter examples to the {\em original} question of Kac | namely in the planar case over $\mathbb{R}$ | arose via that bridge.

\medskip
\subsection{Different directions}

Almost conjugacy has occurred (independently) in various other parts of Mathematics and Physics. For instance,  
let $\K$ be a finite Galois extension of $\mathbb{Q}$ with Galois group $G = \mathrm{Gal}(\K/\mathbb{Q})$, and let $\K_1$ and
$\K_2$ be subfields of $\K$ corresponding to subgroups $G_1$ and $G_2$ of $G$, respectively. Then the zeta functions of $\K_1$ and $\K_2$ are the same if and only if $G_1$ and $G_2$ are almost conjugate in $G$.
In particular, if $G_1$ and $G_2$ are not conjugate in $G$, then $\K_1$ and $\K_2$ are {\em not}
isomorphic while having {\em the same} zeta function. 

The theme of isospectrality of surfaces has occurred, in one of its many guises, and often independently, in (but not restricted to) the following list.

\begin{itemize}
\item
{\em Acoustics}.\quad
\item
{\em Quantum Chaos}.\quad
\item
{\em Nonlinear Dynamics}.\quad
\item
{\em Topology of Manifolds}.\quad
\item
{\em Combinatorial Theory}.\quad
\item
{\em Number Theory}.\quad
\item
{\em Group Theory}.\quad
\item
{\bf $\ldots$}
\end{itemize}

\medskip
\subsection{The EC condition as a setting for construction and classification}

Call a group triple $(G,H,K)$ an {\em EC-triple} if each element of $H$ is conjugate to some element in $H$ (and vice versa). Such triples generalize AC-triples, and $H$ and $K$ are said to be {\em elementwise conjugate} in $G$.

In the same vein as Sunada's Theorem \ref{sunada}, one can obtain the following result.
Let $\M$ be a closed hyperbolic $n$-manifold, $G$ a group, and $H$ and $K$ elementwise
conjugate subgroups of $G$. Then if $\pi_1(\M)$ admits a homomorphism onto $G$, the covers $\M_H$
and $\M_K$ associated to the pullback subgroups of $H$ and $K$ are length equivalent.

For our purposes it is more convenient to first study the EC-property,  since 
it is more flexible in several ways, e.g., for manupilation in induction arguments. As AC implies EC, any classification result on 
the broader class of EC-triples $(G,H,K)$ will yield direct results for the AC-triples. 

But there is more. In some sense, the EC-property reflects a combinatorial property of an underlying incidence geometry (which I have described in this text), which is a crucial tool in the whole theory that is set up in this paper (and we will come back to this geometry in a near future).

\medskip
\subsection{The present paper}

We start this paper with a detailed explanation of the connection between the notions of ``transplantability'' and ``almost conjugacy,'' which are crucial in the construction theory of planar isospectral billiards. We then discuss the related notions {\em AC-triple} (or {\em Gassman-Sunada triple}) and {\em EC-triple} (which naturally generalizes AC-triples), in the context of isospectral drumheads (in any dimension).  Next, we explore some ``naive'' constructions of EC-triples, starting from one given EC-triple, through direct products and ``adding kernels.'' Then we will
introduce four properties | being FF, MAX, PAIR, INV | inspired by physical properties which ``irreducible'' drumheads should have on the one hand, and which planar examples {\em have}, on the other.

We will introduce a  set of extremely general construction procedures of EC-triples  with FF and MAX (starting from a given one), which due to MAX, is connected to the O'Nan-Scott Theorem of finite Group Theory, and which involves finite simple groups. 
Contrary to the fact that most known examples of isospectral manifolds arise from solvable groups (see for instance de Smit and Lenstra \cite{dSL}) | in particular, from two-step nilpotent $(p)$-groups
such as the finite Heisenberg group (see for instance Guralnick \cite{Gur}) | our examples live at the complete other end of the spectrum: they essentially arise from finite simple groups!

\medskip
\begin{center}
\texttt{simple groups}\ +\ \texttt{data}\quad $\longrightarrow$\quad \texttt{irred.}\ \texttt{length equiv.}\ \texttt{manifolds}
\end{center}

\medskip
In a final stage, we will classify EC-triples, by starting from EC-triples with FF, MAX, PAIR and INV, and we will show that, indeed, any such example is constructed through one of our procedures. The property MAX allows us to use the O'Nan-Scott Theorem to some extent.
In particular, all the known planar counter examples to Kac's question ``Can one hear the shape of a drum'' eventually arise, after careful analysis of the properties beyond O'Nan-Scott.

\medskip
\begin{center}
 \texttt{irred.}\ \texttt{length equiv.}\ \texttt{manifolds} \quad $\overset{\texttt{O'Nan-Scott Theory}}{\longrightarrow}$\quad \texttt{simple groups}
\end{center}

\medskip
In such a way, we will obtain a deep connection between irreducible (possibly higher-dimensional) drumheads and finite simple groups that offers a new direction in the theory.

\newpage
\section{Tiles and transplantation}

All known pairs of planar isospectral billiards can be described through Gassman-Sunada group triples which have a very specific 
form. In this section, we analyze this situation. We also clarify the notion of ``transplantation'' in its connection with AC, a subject which seems to be obscured at some places in the literature.

\subsection{Involution graphs}

All known isospectral billiards can be obtained by {\em unfolding polygonal-shaped tiles} | see \cite{GirTha}, but essentially one can only consider triangles. The way the tiles are unfolded can be specified by three permutation $(N\times N)$-matrices $M^{(\mu)}$, $1 \leq \mu \leq 3$ and $N \in \mathbb{N}^{\times}$, associated with the three sides of the triangle: 

\begin{itemize}
\item
$M^{(\mu)}_{ij}  =  1$ if tiles $i$ and $j$ are glued by their side $\mu$;
\item
$M^{(\mu)}_{ii}  =  1$ if the side $\mu$ of tile $i$ is on the boundary of the billiard, and 
\item
$0$ otherwise. 
\end{itemize}
The number of tiles is, of course, $N$.

One can encrypt the action of the $M^{(\mu)}$ ``in'' a graph with colored edges: 
each copy of the base tile is associated with a vertex, and vertices $i$ and $j$, $i \ne j$, are joined by an edge of color $\mu$ if and only if $M^{(\mu)}_{ij}  =  1$. 
In the same way, in the second member of the pair, the tiles are unfolded according to permutation matrices $N^{(\mu)}$, $1 \leq \mu \leq 3$. 
We call such a colored graph an {\em involution graph} (for reasons to be made clear later on).

\br\label{tree}{\rm
As motivated in \cite{GirTha}, we assume the hypothesis that these involutions graphs are {\em trees}.
}
\er

\bigskip
\begin{figure}[h!]
  \centering
    \includegraphics[width=\textwidth]{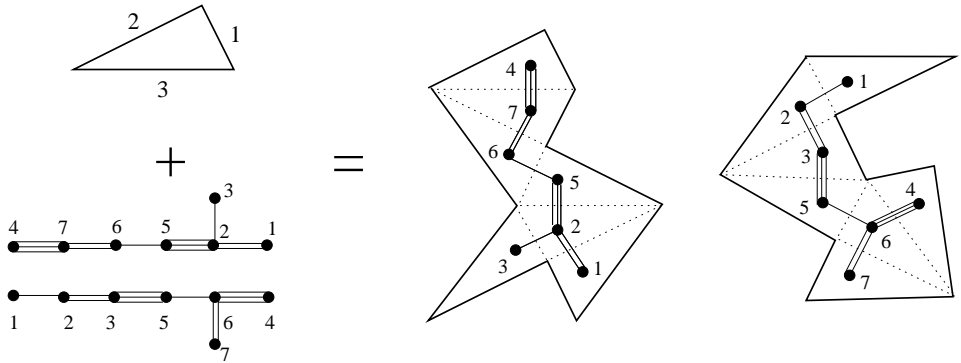}
      \caption{Tiling a triangle through an involution graph; here, the colors are represented by multiple edges.}
\end{figure}

\bigskip
\br{\rm
The consideration of the colored graph as such defined was first made by Y.~Okada and A.~Shudo in \cite{OS}.}\\
\er

\subsection{Transplantability}

Two ``polygonally tiled'' billiards are said to be {\em transplantable} if there exists an invertible matrix $T$ | the {\em transplantation matrix} | such that 
\begin{equation}
 TM^{(\mu)}  =  N^{(\mu)}T \ \ \mbox{for all}\ \ \mu. 
 \end{equation}

Define $G_M := \langle M^{(\mu)}\ \vert\ \mu = 1,2,3 \rangle$ (with ordinary matrix multiplication) and $G_N := \langle N^{(\mu)}\ \vert\ \mu = 1,2,3 \rangle$; then $G_M$ and $G_N$ are subgroups of $\GL_N(\C)$, and the existence of the invertible matrix $T$ {\em precisely} expresses the fact that $G_N$ and $G_M$ are equivalent (= isomorphic) $\GL_N(\C)$-linear representations of the same group. 

\br{\rm
We note that the $M^{(\mu)}$-elements and $N^{(\nu)}$-elements are involutions.
}
\er

On the other hand, since the $M^{(i)}$ and $N^{(j)}$ are permutation matrices, $G_M$ and $G_N$ are subgroups of the symmetric group $\texttt{S}_N$ | that is, they both come with an action on $\Upsilon = \{1,2,\ldots,N\}$. Now the aforementioned linear representations are equivalent if and only if they have the same permutation character | i.e., there exists an isomorphism
\begin{equation}
\chi: G_M \longrightarrow G_N
 \end{equation}
such that 
\begin{equation}
\#\mathrm{Fix}(g) = \#\mathrm{Fix}(g^{\chi})\ \ \forall\ \ g \in G_M.
\end{equation}

Note that the permutation groups $(G_M,\Upsilon)$ and $(G_N,\Upsilon)$ are respectively similar to the actions given by left translation on the left coset spaces $(G_M,G_M/U)$ and $(G_N,G_N/V)$, where $U = {(G_M)}_x$ is an arbitrary point stabilizer in $(G_M,\Upsilon)$, and $V = {(G_N)}_y$ an arbitrary point stabilizer in $(G_N,\Upsilon)$. After applying $T$ (by conjugation), we obtain permutation representations $(G_M,G_M/U)$ and $(G_M,G_M/V^T)$ with 
the same permutation character. \\

Conversely, let $(G,X)$ and $(G,Y)$ be two faithful transitive permutation groups of the same degree. As we have seen, we can identify 
both actions similarly with left translation actions $(G,G/U)$ and $(G,G/V)$ with $U, V$ subgroups of $G$, and $[G : U] = [G : V]$. Now $(G,X)$ and 
$(G,Y)$ have the same permutation character if and only if AC is satisfied for $(G,U,V)$. 

Let $k \in \{ \mathbb{Q},\mathbb{R},\C \}$, and let $V := k[G/H]$ be the $k$-vector space defined by the formal sums over $k$ with base the elements of the left coset space $G/H$, $H$ any subgroup of $G$. Then $G$ acts naturally on $V$ by left translation. If $B$ is the (standard) base 
of $V$ corresponding to $G/H$, it is clear that $G$ naturally defines a subgroup $G_H$ of $\GL_{\ell}(V)$, with $\ell = [G : H]$ (by its action on $B$ and linear extension). 
In this way, we have a faithful $k$-linear representation
\begin{equation}
\rho_H:\ \ G \ \ \hookrightarrow\ \ \GL_{\ell}(V).
\end{equation} 

Now AC is satisfied for $(G,U,V)$ if and only if $\rho_U$ and $\rho_V$ are equivalent, that is, if and only if there exists a $T \in \GL_{\ell}(V)$ such that for all $g \in G$,
\begin{equation}
\rho_U(g)\circ T = T \circ \rho_V(g).
\end{equation}
(Here, $k[G/U]$ and $k[G/V]$ are naturally identified.)

\medskip
\bt[Transplantability]
Two faithful permutation representations $(G,G/U)$ and $(G,G/V)$ of the same degree (with $U$ and $V$ subgroups of the group $G$) have the same permutation character if and only if AC is satisfied for $(G,U,V)$ if and only if the faithful $k$-linear representations $\rho_U$ and $\rho_V$ are $k$-linearly equivalent  if and only there exists a $T \in \GL_{\ell}(V)$ such that for all $g \in G$, 
$\rho_U(g)\circ T = T \circ \rho_V(g)$.
\et

In other words: almost conjugacy and transplantability express {\em the same property}.

Abstractly, we will  sum up some equivalent properties in the next section.

\subsection{Connection with shape}

If the matrix $T$ is a permutation matrix, the two domains would just have the same shape (i.e., are isometric), and vice versa.

\subsection{Analysis of the known examples}

All {\em known} planar counter examples $(\mD_U,\mD_V)$ to Kac's question are transplantable, i.e., come from Gassman-Sunada data $(G,U,V)$; see e.g. \cite{GirTha}. Much more can be said about the triples.

\begin{itemize}
\item[(i)]
In all
known examples, the group $G$ is isomorphic to the classical simple group $\PSL_n(q)$, where $(n,q) \in \{(3,2),(3,3),$ $(4,2),(3,4)\}$, cf. \cite{Conw}. 
\item[(ii)]
In all these examples, $U$ and $V$ are {\em maximal subgroups} of $G$.
\item[(iii)]
Also, the 
$T$-elements are involutions which play the role of a duality of a naturally associated projective space $\PG(n - 1,q)$, on which the $M^{(\mu)}$-elements, the
$N^{(\nu)}$-elements and $\PSL_n(q)$ naturally act as automorphisms.
\item[(iv)]
By combined work of O. Giraud and K. Thas, it follows that, conversely, all transplantable pairs arising from the special linear group, and satisfying Equation (\ref{Fixeq}) below, are the known ones.
This was observed in \cite{Giraud,KTI,KTII,KTIII}.
\end{itemize}

We also mention the following result.

\begin{lemma}[Y. Okada and A. Shudo \cite{OS}]
\label{Shudo}
All isospectral transplantable drums, that unfold an $r$-gon $N$ times, are known if $N \leq 13$.
\end{lemma}

\subsection*{Fixed points equation}

The involution graph of $\mD_U$ (respectively $\mD_V$) which was earlier defined is  the so-called (undirected) ``Schreier coset graph'' of $G$ with respect to $U$ (respectively $V$)
relative to $\{M^{(\mu)}\}$ (respectively $\{ N^{(\nu)}\}$) | which is a generating set of involutions of $G$ | by definition.

As the Schreier graphs under consideration are trees (cf. Remark \ref{tree}), the following identity holds:
\begin{equation}\label{Fixeq}
(r-2)\Lambda =\sum_{j=1}^{r} \mbox{ Fix }(\chi^{(j)})-2,
\end{equation}
where Fix$(\chi^{(j)})$ is the number of fixed points of $\chi^{(j)}$ in the permutation group $(G_{\chi},G_{\chi}/S)$, with $\chi = M, N$ and $S = U, V$ respectively; see \cite{GirTha}. Here, $\Lambda$ is the number of tiles, and $r$ is the gonality of the base polygon, which may be taken to equal $3$.

\section{Elementwise and almost conjugacy}

The following properties are easily shown to be equivalent. Let $H, H'$ be subgroups of a finite group $G$.
\begin{itemize}
\item[(a)]
$H$ and $H'$ are almost conjugate (AC) (that is, $(G,H',H')$ is a Gassmann-Sunada triple).
\item[(b)]
There exists a bijection $B: H \longrightarrow H'$ such that for each $h \in H$, $B(h)$ is conjugate to $h$ in $G$.
\item[(c)]
For all $g \in G$, we have that $\vert g^G \cap H \vert = \vert g^G \cap H' \vert$.
\end{itemize}

\begin{theorem}[AC and Isospectrality]
\label{ACIso}
Let $\M$ be a closed hyperbolic $n$-manifold, $G$ a group, and $H$ and $K$ almost
conjugate subgroups of $G$. Then if $\pi_1(\M)$ admits a homomorphism onto $G$, the covers $\M_H$
and $\M_K$ associated to the pullback subgroups of $H$ and $K$ are isospectral.
\end{theorem}

We say that subgroups $H$ and $H'$ of a finite group $G$ are {\em elementwise conjugate} if for each $g \in G$ we have that
\begin{equation}
g^G \cap H \ne \emptyset \ \   \Longleftrightarrow\ \      g^G \cap H' \ne \emptyset
\end{equation}

The following properties are  equivalent.
\begin{itemize}
\item[(i)]
$H$ and $H'$ are almost elementwise conjugate (EC).
\item[(ii)]
Each element of $H'$ is conjugate to an element of $H$ and each element of $H$ is conjugate to some element of $H'$.
\end{itemize}

For $\M$ a compact Riemannian manifold, let $\mathcal{L}(\M)$ be the set of lengts of closed geodesics of $\M$ (that is, $\mathcal{L}(\M)$ is the {\em length spectrum} without multiplicities). Then two such manifolds $\M$ and $\widetilde{\M}$ are {\em length equivalent} if $\mathcal{L}(\M) = \mathcal{L}(\widetilde{\M})$.

Expressed in terms of isospectrality-related properties, we have the following version of Theorem \ref{ACIso}.

\begin{theorem}[EC and Length Equivalence]
\label{ECLE}
Let $\M$ be a closed hyperbolic $n$-manifold, $G$ a group, and $H$ and $K$ elementwise
conjugate subgroups of $G$. Then if $\pi_1(\M)$ admits a homomorphism onto $G$, the covers $\M_H$
and $\M_K$ associated to the pullback subgroups of $H$ and $K$ are length equivalent.
\end{theorem}

The following is obvious.
\begin{observation}
\label{EXAC}
Let $(G,H,H')$ be a Gassmann-Sunada triple. Then $H$ and $H'$ are elementwise conjugate in $G$.
\end{observation}

So in general we have
\begin{equation}
\mathrm{AC} \ \ \longrightarrow\ \ \mathrm{EC}.
\end{equation}

The converse is not necessarily true | one cannot even conclude that $H$ and $H'$ have the same order!
For example, let $A_n$ be the alternating group on $n$ letters ($n \geq 3$), and consider commuting involutions (transpositions) $\alpha$ and $\beta$.
Then $H = \{ \id,\alpha \}$ and $H' = \{ \id,\alpha,\beta,\alpha\beta \}$ are EC (and certainly not AC).

%\begin{remark}
%{\rm 
%I have seen in some places that 
%}
%\end{remark}

Still, for our purposes it is more convenient to first study the EC-property,  since 
it is more flexible in several ways (e.g. for manupilation in induction arguments). As AC implies EC, any classification result on 
the broader class of EC-triples $(G,H,H')$ will yield direct results for the AC-triples.

\section{Construction of isospectral drums}

In this section, we want to construct isospectral  drums (not necessarily planar) arising from  Gassman-Sunada group triples $(G,H,H')$. 
We will define a class $\mD$ of drums which contains all presently known examples of isospectral {\em planar} drums. Our way to construct starts from a known example, and makes an infinite class of examples out of it. 

We start very generally, with generic constructions which give rise to many EC-triples. We always start from some Gassman-Sunada triple, or more generally with an EC-triple.
We then fine-tune our construction in such  a way that
the drums which come out satisfy some strong properties which are suggested by the physics of the game (cf. the next subsection). On the other hand, we want every drum/triple which is constructed to inherit certain properties from the initial example. 
This also means 
that an initial triple $(G,H,H')$ has to satisfy a number of properties in order to make it eligible for our method to apply. Our model example will be the projective linear groups $\PGL_n(q)$ over a finite field $\F_q$. 

So the idea is this: in each step of our theoretic construction, the construction becomes more and more general, and produces a wealth of (new) examples of triples/drums. But in order to apply it,
in each step, the needed base triples become more specialized (but at this point only slightly). After having established the construction method  the connection with simple groups will be revealed  in the next section,.

\subsection{Properties FF, MAX, INV and PAIR}

We will give extra 
attention to the following properties, which are partly motivated by the list of known isospectral planar drums.\\

Recall that a subgroup $U, \cdot$ of a group $V, \cdot$ is {\em normal} if for every $(u,v) \in U \times V$, we have that 
\begin{equation}
v^{-1}\cdot u\cdot v \in U.
\end{equation}
(Below, we will omit the ``$\cdot$-notation'' if the group operation is clear.)

\begin{itemize}
\item[FF]
The FF-property expresses the fact that the actions $G \curvearrowright G/H$ and $G \curvearrowright G/H'$ are {\em faithful}. This means that there are no nontrivial normal
subgroups of $G$ contained in $H$ and $H'$. Note that if $N$ is a normal subgroup of $G$ inside $H$, $N$ is also contained in $H'$, since each element of $H'$ is conjugate to some element in $H$, and $N$ is self-conjugate.
This assumption is natural, since 
dividing out normal subgroups yields the same drums. (``The drums do not see these normal subgroups.'')\\
\end{itemize} 

\bpr[FF]
By its mere definition, any example of isospectral drums coming from involution graphs \`{a} la Okada-Shudo is FF.\eopr \\
\epr

\begin{itemize}
\item[MAX]
MAX expresses the fact that both $H$ and $H'$ are {\em maximal subgroups} of $G$, i.e., if $A$ is a subgroup of $G$ containing $H$ (respectively $H'$), then either $A = G$ or $A = H$ (respectively $A = G$ or $A = H'$). This is an irreducibility assumption, as in some way MAX-examples represent irreducible examples; if for instance $H$ would not be maximal, it might happen that there are proper subgroups $A, B$ of $G$
for which $H \leq A \ne H$ and $H' \leq B \ne H'$ such that $(G,A,B)$ is also Gassman-Sunada/EC. And then physically, the drums coming from the data $(G,H,H')$ consist of the smaller drums coming from the data $(G,A,B)$. \\
\end{itemize}

\bpr[MAX]
Asking that $H$ and $H'$ are maximal means that the drums are not built from smaller ones.\eopr \\
\epr

\begin{itemize}
 \item[INV]
By Involution Formula, we known that transplantable {\em planar} drums satisfy equation (\ref{Fixeq}). The existence of  involutions satisfying (\ref{Fixeq}) is INV.
\item[PAIR]
In the known planar examples, there is an outer automorphism of $G$ which maps $H$ to $H'$ | in those cases, $G$ is isomorphic to a general special projective linear group over some finite field, and $H$ is the stabilizer of some point (in the action on the associated projective space), and $H'$ the stabilizer of some hyperplane, cf. \S\S \ref{model1} and \ref{model2}. The outer automorphism is associated to a duality of the projective space. PAIR expresses this fact: PAIR is satisfied for $(G,H,H')$ if there is some $\sigma \in \mathrm{Out}(G)^\times$ such that $H^\sigma = H'$
and $\sigma^2 \in \mathrm{Inn}(H)$.  In fact, in most cases in which we will need PAIR, much more general versions suffice. As we will see, often it will be enough to ask that $\vert H \vert = \vert H'\vert$.
\end{itemize}

%\medskip
%\subsection{Passing to subgroups}

\medskip
\subsection{Adding kernels}

Let $(S,L,L')$ be a Gassman-Sunada triple which is FF. Let $K$ be any group, and consider $(S \times K,L \times K,L' \times K)$. Here $A \times B$, with $A$ and $B$ denoting groups, is the {\em direct product} of $A$ and $B$. It consists of all elements of the Cartesian product of $A$ and $B$, and the group law is
\begin{equation}
(a,b)\cdot(a'b') := (aa',bb').
\end{equation}
Note that conjugation works as follows: 
\begin{equation}
(a,b)^{(u,v)} = (u,v)^{-1}(a,b)(u,v) = (u^{-1}au,v^{-1}bv) = (a^u,b^v).
\end{equation}

As $(S,L,L')$ is Gassman-Sunada, there is a bijection $B: L \longrightarrow L'$ such that for $s \in L$, $B(s)$ is conjugate to $s$ (in $S$).
Now define the bijection $\overline{B}$ as
\begin{equation}
\overline{B}: L \times K \longrightarrow L' \times K: (s,k) \longrightarrow (B(s),k).
\end{equation}
We know that $B(s) = s^g$ for some $g \in S$; it follows that $(B(s),k) = (s,k)^{(g,1)}$, so that $(S \times K,L \times K,L' \times K)$ clearly 
also is Gassman-Sunada.

But of course, no new information is added, as  $(S \times K,L \times K,L' \times K)$ is not FF | $\{ \id\} \times K$ obviously is a normal subgroup of $S$
and it is contained in both $L \times K, L' \times K$.

\begin{remark}[EC-Triples]
{\rm 
Exactly the same procedure also works when we start with EC-triples. And FF remains valid throughout the construction. }
\end{remark}

\begin{remark}[Arbitrary products]{\rm 
Note that $K$ is arbitrary, so can be replaced by any finite product of finte groups.}
\end{remark}

We will change this procedure slightly, so that we can sustain the FF-property when passing to ``bigger examples.''\\

\subsection{Passing to products}

Start again from data $(S,L,L')$, and consider $(S \times S,L\times L,L'\times L')$. Suppose that $(S,L,L')$ is FF.

Let $B$ be as above, and define $\overline{B}$ as follows:
\begin{equation}
\overline{B}: L \times L \longrightarrow L' \times L': (s,t) \longrightarrow (B(s),B(t))
\end{equation}
We know that there are $g, h$ in $S$ such that $B(s) = s^g$ and $B(t) = t^h$, so $(B(s),B(t)) = (s,t)^{(u,v)}$, and $\overline{B}$ obviously is a bijection. 
It follows that $(S \times S,L\times L,L'\times L')$ is Gassman-Sunada.

The same can be concluded when starting from an EC-triple $(S,L,L')$ (i.e., $(S \times S,L \times L,L' \times L')$ is then EC).

We need to check FF, so suppose that $N$ is a normal subgroup of $S \times S$ contained in $L \times L$ (and so also in $L' \times L'$). Let $\pi_1$ be the projection
of $N$ onto the first component in $S \times S$ (that is, $\pi_1(n,m) = n$ for $(n,m)$ in $S \times S$), and let $\pi_2$ be the projection on the second component. 
Suppose one of the images, say $\pi_1(N)$ is different from $\{\id\}$ and $L$; then $\pi_1(N)$ obviously is a normal subgroup of $S$ inside $L$, contradicting the fact that $(S,L,L')$ is FF. If this property does {\em not} hold for both $\pi_1$ and $\pi_2$, then $\pi_i(N)$ must be either $\{\id\}$ or $L$, and for at least one component, say the second, we have $\pi_2(N) = L$. And then $L$ is a normal subgroup of $S$. 

If $(S,L,L')$ is Gassman-Sunada, we have a contradiction as $L$ and $L'$ are AC, so that $L \unlhd S$ would imply $L = L'$.

Now let $(S,L,L')$ be an EC-triple, and let $L$ be a normal subgroup of $S$. Let $l' \in L'\setminus L$; as $l'$ is conjugate to some element $l$ in $L$, we obtain a contradiction as $l^S \subseteq L$. As $N \leq (L \times L) \cap (L' \times L')$, $\pi_2(N) = L \leq L'$.
We thus conclude that $L = L'$, contradiction.

%\begin{remark}[EC-Triples]
%{\rm 
%Again the same procedure also works when we start with EC-triples. And FF also remains valid throughout the construction.}
%\end{remark}

\medskip
\bpr[Inheritance of FF]
\label{sgr}
Let $(S,L,L')$ be an EC-triple (so that as a particular case we have Gassman-Sunada triples); then
 our construction is guaranteed to give us an example with FF. \eopr \\
\epr

And clearly, instead of considering a product with two components, we can do exactly the same reasoning on {\em any finite number of components}.

\medskip
\bo
If $(S,L,L')$ is Gassman-Sunada with FF, or more generally $(S,L,L')$ is an EC-triple with FF, and $k \ne 0$ is any positive integer, then $(S^{\times k},L^{\times k},{L'}^{\times k})$ is also Gassman-Sunada with FF.\eop \\
\eo

Here, $A^{\times k}$ denotes the direct product $A \times \cdots \times A$ with $k$ factors.

But there is a catch: even if we would assume that $L$ and $L'$ are maximal subgroups of $S$, $L \times L$ and $L' \times L'$ certainly are not, as 
for example $L \times L \leq L \times S$. So MAX is not inherited if we pass to products. 

We need a more subtle approach.\\

\subsection{Products with a twist}

Now consider an EC-triple $(S,L,L')$, with FF (MAX is not needed {\em for now}).  We pass to $(S \times S,L \times L,L'\times L')$ as above, which is also FF by Principle \ref{sgr}. Now consider $\mathrm{S}_2 = \langle \sigma \rangle =: B$, the symmetric group acting on $\{1,2\}$  ($\sigma = \sigma^{-1}$). If $(s_1,s_2) \in S \times S$, we define $\sigma((s_1,s_2))$ as 
\begin{equation}
(s_{\sigma^{-1}(1)},s_{\sigma^{-1}(2)}) = (s_2,s_1).
\end{equation}
 
Put $A := S \times S$, and define a group $A \rtimes B$ as follows. Its elements are just of the form $(a,b)$ with $a \in A$ and 
$b \in B$, but the group operation is given by
\begin{equation}
(a,b)\cdot(a',b') := (ab(a'),bb').
\end{equation}
By definition, $A \rtimes B$ is the {\em wreath product} $S \wr B$ of $S$ by $B$. 
Note that $A$ and $B$ naturally are (isomorphic to) subgroups of $A \rtimes B$; note also that $A$ is (isomorphic to) a {\em normal} subgroup of $A \rtimes B$.\\

Let $(l,r), (n,h)$ be elements of $A \rtimes B$; then 
\begin{equation}
\label{eqeq}
(n,h)^{(l,r)} = (l,r)^{-1}(n,h)(l,r) = (r^{-1}(l^{-1}),r^{-1})(n,h)(l,r) = (r^{-1}(l^{-1}nh(l)),h^r).
\end{equation}
In particular, if $r = 1$, then we obtain that $(n,h)^{(l,r)} = (l^{-1}nh(l),h)$.\\

Define subgroups $H$ and $H'$ of $G := A \rtimes B$ by $H := (L \times L) \rtimes B$ and $H' := (L' \times L') \rtimes B$.
It is clear that $L \times L$ and $L' \times L'$ can be naturally seen as subgroups of index two of respectively $H$ and $H'$.

Now let $\overline{B}: L \times L \longrightarrow L' \times L'$ be as before, and define $\widetilde{B}$ as follows:
\begin{equation}
\widetilde{B}: H \longrightarrow H': (n,b) \longrightarrow (\overline{B}(n),b).
\end{equation}
Here, $n \in L \times L$ and $b \in \langle \sigma \rangle = B$. The fact that $\overline{B}$ is a bijection from $L \times L$ to $L' \times L'$
readily implies that $\widetilde{B}$ is a bijection between $H$ and $H'$. \\

To show that $H$ and $H'$ are EC in $G$, we need to show that for all $h' \in H'$, we can find $u \in G$ and $h \in H$ such that 
$h^u = h'$. By the above, this property is certainly true for the subgroup $L' \times L'$ of $H'$ (each of its elements is conjugate in $G$ to 
some element of $L \times L$). So it is sufficient only to consider elements of type 
$(n,\sigma)$ in $H'$ (where $n = (a,b) \in L' \times L'$). The symmetrical argument then yields that the obtained triple is EC.
So we seek for $(l,\delta) \in G$ (with $l = (l_1,l_2) \in S \times S$ and $\delta \in \{\id,\sigma\}$) and $(r,\sigma)$
(with $r \in L \times L$) such that
\begin{equation}
(r,\sigma)^{(l,\delta)} = (n,\sigma).
\end{equation}

Plugging this into equation (\ref{eqeq}), we obtain the following equations, the first being for $\delta = \sigma$, the second for 
$\delta = \id$, and the third being a trivial equation:

\begin{equation}\left\{
\begin{array}{ccccc}
r &= & l\sigma(n)\sigma(l)^{-1} &= & (l_1bl_2^{-1},l_2al_1^{-1})\\
r &= & ln\sigma(l)^{-1} & = & (l_1al_2^{-1},l_2bl_1^{-1})\\
\sigma &= &\sigma.\\
\end{array}\right.
\end{equation}

\medskip
\bpr[Reduction to the base group $S$]
Since $(a,b) \in L' \times L'$ is essentially arbitrary, both equations boil down to one: we need $(l_1,l_2) \in S \times S$ such that
$(l_1al_2^{-1},l_2bl_1^{-1}) \in L \times L$. In fact, since $l_2bl_1^{-1} \in S$ if and only if $l_1b^{-1}l_2^{-1} \in S$, and since $b$
is arbitrary,  we need $(l_1,l_2) \in S \times S$ such that
$(l_1al_2^{-1},l_1bl_2^{-1}) \in L \times L$.
\eopr \\
\epr

\br[EC and its geometry]{\rm
At this point, one has to start specifying the initial triple $(S,L,L')$ in order to play the game.
We will introduce a fundamental example in the next
subsection. Much seems to depend on detailed knowledge of the natural (geometric) modules on which the groups act. However, we will
surprisingly observe that this is not true, and that the geometric property we need naturally comes from the EC-property (in fact, it is even {\em equivalent}). This is another
reason why the class of EC-triples is a more natural candidate to study for our purposes than the class of AC-triples.
}
\er

\medskip
\subsection{Model example: the general linear groups}
\label{model1}

Consider the simple group $\PGL_n(\F)$ (over any commutative field $\F$, $n \in \mathbb{N}^\times \setminus \{1\}$), and put $S = \PGL_n(\F)$. As we have seen earlier, $S$ acts naturally on the linear subspaces of the projective space $\mathbf{PG}(n - 1,\F)$.
Let $L := S_x$ be the stabilizer in $S$ of a point $x$ in $\mathbf{PG}(n - 1,\F)$, and let $L' = S_{\Pi}$ be the stabilizer in $S$ of a 
hyperplane $\Pi$ in $\mathbf{PG}(n - 1,\F)$. It is well known and easy to see that $(S,S_x,S_{\Pi})$ has the EC-property. (If $n = 2$, $S_x$ and $S_{\Pi}$ are conjugate, though.)

\br{\rm
The choice of $(x,\Pi)$ is not important, as the reader will see.
}
\er

So $a, b$ both stabilize $\Pi$, and hence $ab^{-1}$ has the same property. 
A well-known property for projective spaces tells us that $ab^{-1}$ therefore must fix some point $y$ of the space. 
This means that
\begin{equation}
a(y) = b(y) = y'
\end{equation}
for some point $y'$. \\

\medskip
\begin{center}
\begin{tikzpicture}
\fill [ultra thick,blue] (3,2) circle (2pt);
\fill [ultra thick,blue] (7,2) circle (2pt);
\fill [ultra thick,blue] (5,0) circle (2pt);
\draw (2,2) node[draw=none,fill=none] {point $y$};
\draw (5,3.3) node[draw=none,fill=none] {$a$};
\draw (5,0.7) node[draw=none,fill=none] {$b$};
\draw (3.5,0.5) node[draw=none,fill=none] {$l_1^{-1}$};
\draw (6.5,0.5) node[draw=none,fill=none] {$l_2$};
    \draw[->] [ultra thick,blue] (3,2) to[out=45,in=135]  (7,2);
        \draw[<-] [ultra thick,blue] (3,2) to[out=-45,in=-135]  (7,2);
        \draw (8,2) node[draw=none,fill=none] {point $y'$};
          \draw[->] [ultra thick,blue] (3,2) to[out=-45,in=-135]  (5,0);
        \draw[->] [ultra thick,blue] (5,0) to[out=-45,in=-135]  (7,2);
        \draw (5,-0.75) node[draw=none,fill=none] {point $x$};
       % \draw [decoration={
          %   markings,
             %mark=at position 0.4 with \arrow{>}}
       %,postaction=decorate]
       %(A) to (B);
\end{tikzpicture}
\end{center}
%(5,3) to[out=-180+115,in=10]

Let $l_1^{-1}$ map $y$ to $x$ (note that such elements certainly exist!), and let $l_2$ map $x$ to $y'$ (same remark). A standard property
of permutation groups now tells us that all elements in $S$ that map $y$ to $y'$ are in the set
\begin{equation}
l_1^{-1}S_xl_2,
\end{equation}
and hence $(l_1al_2^{-1},l_1bl_2^{-1}) \in S_x \times S_x = L \times L$. 

\bp[Linear groups and $2$-cycles]
Put 
\begin{equation}
(S,L,L') = (\PGL_n(\F),{\PGL_n(\F)}_x,{\PGL_n(\F)}_{\Pi}) 
\end{equation}
with $(x,\Pi)$ any point-hyperplane pair. Then 
$(S \wr S_2,L \wr S_2,L' \wr S_2)$ is an EC-triple. 
\eop 
\ep

We will observe further on how this construction behaves with respect to the properties we introduced in the beginning of this section.
First, we want to replace the involutive action by a $3$-cycle action. This will lead us to the general construction, which will handle 
wreath products of simple groups {\em with any transitive subgroup of the symmetric group}.

\medskip
\subsection{$3$-Products with a $3$-cycle}

Consider an EC-triple $(S,L,L')$, with FF (MAX is again not asked).  We pass to $(S \times S \times S,L \times L \times L,L'\times L' \times L')$ as above, which is also FF by Principle \ref{sgr}. Now consider any $\gamma \in \mathrm{S}_3$, the symmetric group acting on $\{1,2,3\}$. 
If $(s_1,s_2,s_3) \in S \times S \times S$, we define $\gamma((s_1,s_2,s_3))$ as 
\begin{equation}
(s_{\gamma^{-1}(1)},s_{\gamma^{-1}(2)},s_{\gamma^{-1}(3)}).
\end{equation}
 
Put $A := S \times S \times S$, and define a group $S \wr B$ as in the previous subsection (see also the next remark), where $B \leq \mathrm{Sym}(3)$ is supposed to be transitive. 

\begin{remark}[Wreath products]{\rm
Let $G$ be any group, and let $P$ be a subgroup of the symmetric group $\mathrm{Sym}(\Omega)$ acting on the finite set $\Omega \ne \emptyset$. 
The {\em wreath product} $G \wr P$ of $G$ by $P$ is defined as follows.
Its elements are of the form $(g,p)$ with $g \in G^\Omega$ (which is the Cartesian product of $\vert \Omega\vert$ copies of $G$, indexed
by $\Omega$) and 
$p \in P$, while the group operation is given as before by
\begin{equation}
(g,p)\cdot(g',p') := (gp(g'),pp').
\end{equation}
Here, with $g' = {(g_\omega)}_{\omega \in \Omega}$, 
\begin{equation}
p(g') := {(g_{p^{-1}(\omega)})}_{\omega \in \Omega}.
\end{equation}
One  notes that $G$ and $P$ are naturally isomorphic to  subgroups of $G \wr P$, and that $G$  becomes a normal subgroup under this
natural identification.\\
}
\end{remark}

To show that $H$ and $H'$ are EC in $G$, we need to show that for all $h' \in H'$, we can find $u \in G$ and $h \in H$ such that 
$h^u = h'$. By the above, this property is certainly true for the subgroup $L' \times L' \times L'$ of $H'$ (each of its elements is conjugate in $G$ to 
some element of $L \times L \times L$, even in a bijective correspondence). So it is sufficient only to consider elements of type 
$(n,\gamma)$ in $H'$ (where $n = (a_1,a_2,a_3) \in L' \times L' \times L'$). So we seek for $(l,\delta) \in G$ (with $l = (l_1,l_2,l_3) \in S \times S \times S$ and $\delta \in B$) and $(r,\gamma')$
(with $r \in L \times L \times L$ and $\gamma' \in B$) such that
\begin{equation}
(r,\gamma')^{(l,\delta)} = (n,\gamma).
\end{equation}

Plugging this into equation (\ref{eqeq}), we obtain the following equations, the first being for $\delta = \sigma$, the second for 
$\delta = \id$. (For general $\delta \in B$, one needs that ${\gamma'}^\delta = \gamma$, but we will show below that the case $\delta = \id$ and 
$\gamma' = \gamma$
already suffices to find solutions.) So below, we let $\gamma' = \gamma$.

\begin{equation}\left\{
\begin{array}{ccccc}
r &= & l\sigma(n)\sigma(l)^{-1} &\\
r &= & ln\sigma(l)^{-1} & = & (l_1a_1l_{\gamma^{-1}(1)}^{-1},l_2a_2l_{\gamma^{-1}(2)}^{-1},l_3a_3l_{\gamma^{-1}(3)}^{-1}).\\
\end{array}\right.
\end{equation}

\medskip
\subsection{Model example: the general linear groups}
\label{model2}

We show again that with $S = \PGL_n(q)$ ($n$, $q$ as above), $L = S_x$ and $L' = S_{\Pi}$ ($x$ and $\Pi$ as above), we can make the desired properties work. We choose an ``arbitrary'' $\gamma \in \mathrm{S}_3$, but we assume it to be a regular $3$-cycle (that is, it generates a sharply transitive subgroup of $\mathrm{S}_3$). The reason is that if it is not of this form, the existence  of the desired $(l_1,l_2,l_3)$ follows from 
it cycle decomposition together with the fact that we already obtained the result for $n = 2$. In any case, the technique we present also works for {\em any} cycle.

Choose for instance (and without loss of generality) the permutation

\begin{equation}\left\{
\begin{array}{ccc}
\gamma(1) &= &3\\
\gamma(2) &= &1\\
\gamma(3) &= &2,\\
\end{array}\right.
\end{equation}
so that we must check one of the conditions 

\begin{equation}\left\{
\begin{array}{ccccc}
r &= & l\gamma(n)\gamma(l)^{-1} & &\\
r &= & ln\gamma(l)^{-1} & = & (l_1a_1l_{2}^{-1},l_2a_2l_{3}^{-1},l_3a_3l_{1}^{-1}).\\
\end{array}\right.
\end{equation}

We solve, as before, the second set of equations.

As $a_1,a_2,a_3$ are elements of $S_{\Pi}$, their inverses are in $S_{\Pi}$ as well, and so $a_2 \circ a_1^{-1} \circ a_3$ is also in $S_{\Pi}$.
As before, we know this element fixes some point, say $y$. Put $a_3(y) =: y_1$ and $a_1 \circ a_3(y) =: y_2$; then $a_2(y_2) = y$.
Now choose elements $l_1, l_2, l_3$ such that:

\begin{equation}\left\{
\begin{array}{ccc}
l_1(x) &= &y_1\\
l_2(x) &= &y_2\\
l_3(x) &= &y,\\
\end{array}\right.
\end{equation}

\medskip
\begin{center}
\begin{tikzpicture}
\fill [ultra thick,blue] (3,2) circle (2pt);
\fill [ultra thick,blue] (7,2) circle (2pt);
\fill [ultra thick,blue] (5,0) circle (2pt);
\fill [ultra thick,blue] (5,3.5) circle (2pt);
\draw (2,2) node[draw=none,fill=none] {point $y$};
\draw (3.5,3.2) node[draw=none,fill=none] {$a_3$};
\draw (6.6,3.2) node[draw=none,fill=none] {$a_1$};
\draw (5,0.7) node[draw=none,fill=none] {$a_2$};
\draw (3.5,0.5) node[draw=none,fill=none] {$l_3$};
\draw (6.5,0.5) node[draw=none,fill=none] {$l_2$};
\draw (4.8,2) node[draw=none,fill=none] {$l_1$};
\draw (5,4.2) node[draw=none,fill=none] {point $y_1$};
    \draw[->] [ultra thick,blue] (3,2) to[out=45,in=135]  (5,3.5);
       \draw[->] [ultra thick,blue] (5,3.5) to[out=45,in=135]  (7,2);
        \draw [ultra thick,blue] (3,2) to[out=-45,in=-135]  (7,2);
        \draw (8,2) node[draw=none,fill=none] {point $y_2$};
          \draw[<-] [ultra thick,blue] (3,2) to[out=-45,in=-135]  (5,0);
            \draw[<-] [ultra thick,blue] (5,3.5) to[out=-45,in=135]  (5,0);

        \draw[->] [ultra thick,blue] (5,0) to[out=-45,in=-135]  (7,2);
        \draw (5,-0.75) node[draw=none,fill=none] {point $x$};
       % \draw [decoration={
          %   markings,
             %mark=at position 0.4 with \arrow{>}}
       %,postaction=decorate]
       %(A) to (B);
\end{tikzpicture}
\end{center}

\medskip
\subsection{The general case: general twists}
\label{gcgt}

Consider an EC-triple $(S,L,L')$, with FF (MAX is again not needed for the construction). Take any $n \in \mathbb{N}^\times$. 
We pass to $(S^{\times n},L^{\times n},{L'}^{\times n})$ as above, which is also FF. Now we take any transitive subgroup $B$
in the symmetric group $\mathrm{S}_n$, acting naturally  on $\{1,\ldots,n\}$. 
If $(s_1,\ldots,s_n) \in S^{\times n}$, we define $\gamma((s_1,\ldots,s_n))$ for any $\gamma \in B$ as 
\begin{equation}
(s_{\gamma^{-1}(1)},\ldots,s_{\gamma^{-1}(n)}).
\end{equation}
 
Put $A := S^{\times n}$, and define $G := S \wr B$, $H := L \wr B$ and $H' :=  L' \wr B$ as before. 

To show that $H$ and $H'$ are EC in $G$, we consider any element $((a_1,\ldots,a_n),\gamma) \in H' = {L'}^{\times n} \rtimes B$, and we show that
it is conjugate to some $(r,\gamma) \in H = L^{\times n} \rtimes B$ through some element of the form $((l_1,\ldots,l_n),\id)$. It is not necessary to 
consider elements with trivial $\gamma$, since they already have the required property.

So we want to solve the set of equations

\begin{equation}\label{eqgamma}
\left\{
\begin{array}{ccccc}
%r &= & l\sigma(n)\sigma(l)^{-1} &= & (l_1bl_2^{-1},l_2al_1^{-1})\\
r &= & ln\gamma(l)^{-1} & = & (l_1a_1l_{\gamma^{-1}(1)}^{-1},\ldots,l_na_nl_{\gamma^{-1}(n)}^{-1}).\\
\end{array}\right.
\end{equation}

\medskip
\subsection{General theory: solving equation (\ref{eqgamma}), and EC-geometry}

Before proceeding, let $(C,D,D')$ be an EC-triple. Define ``points'' as being the elements of the left coset space $\mathcal{X} = \{cD \vert c \in C\}$ and 
``hyperplanes'' as being elements of the left coset space $\mathcal{X}' = \{cD' \vert c \in C\}$. In the left action of $C$ on $\mathcal{X}$, a point 
stabilizer (of the point $gD$) has the form $D^{g^{-1}}$, and in the left action of $C$ on $\mathcal{X}'$, the stabilizer of $hD'$
is ${D'}^{h^{-1}}$. In particular, $D$ is the stabilizer of the point $(\id \cdot) D$ and $D'$ is the stabilizer of the point $(\id \cdot) D'$.
Now by EC, any $a \in D'$ is in some ${D}^{l^{-1}}$, so that $a$ fixes the point $lD$. 

\bpr
EC precisely expresses the geometric property we used in the model example of projective general linear groups to solve the necessary equations for obtaining EC.
\eopr
\epr

Now we return to the setting of \S\S \ref{gcgt}, and let $S, L, L', G, H, H', \gamma$, etc. be as before. We need to solve equation set (\ref{eqgamma}).

By induction on the decomposition cycles of the permutation $\gamma$, it suffices to consider only regular $n$-cycles, although our method 
also applies to non-regular permutations. So from now on, let $\gamma$ be a regular $n$-cycle. As 
$a_1,\ldots,a_n \in S'$, the (composition) product (acting on the left)
\begin{equation}\large
\prod_{i = 1}^na_{\gamma^{-(n - i)}(1)} =: a
\end{equation}
is also in $S'$, and so fixes some point $y$ of the left coset space $\{ sL \vert s \in S \}$ by EC. For each $i = 1,\ldots,n$, we define 
\begin{equation}
\displaystyle\large
a_i(y_i) =: y_{\gamma^{-1}(i)},
\end{equation}
and note that $y_n = y$ (and $\gamma^0 = \id$). We stress the fact that all considered points are elements of the space $\{ sL\ \vert\ s \in S\}$. Put $x$ equal to the point $(\id)\cdot L$, so that $L$ is the point stabilizer  of $x$ in $S$ acting on $\{ sL\ \vert\ s \in S\}$.

\br
{\rm
Note that there could be $i$ such that $y_i = y_{i + 1}$. This phenomenon gives no obstructions. (And this is why the argument works for any type of permutation.)
}
\er

Now for each $i = 1,\ldots,n$, define $l_i$ as being any element of $S$ which sends $x$ to $y_l$ (and note that such elements exist
since $S$ acts transitively on the points of $\{ sL \vert s \in S \}$).\\

\bigskip
\begin{center}
\begin{tikzpicture}
\fill [ultra thick,blue] (2,2) circle (2pt);
\fill [ultra thick,blue] (8,2) circle (2pt);
\fill [ultra thick,blue] (5,0) circle (2pt);
\fill [ultra thick,blue] (5,2) circle (2pt);

\draw (1.2,2) node[draw=none,fill=none] {$y$};
\draw (3.5,3.2) node[draw=none,fill=none] {$a_n$};
\draw (6.6,3.2) node[draw=none,fill=none] {$a_{\gamma^{-1}(n)}$};
%\draw (5,0.7) node[draw=none,fill=none] {$a_2$};
\draw (2.7,0.5) node[draw=none,fill=none] {$l_n$};
\draw (7.8,0.5) node[draw=none,fill=none] {$l_{\gamma^{-2}(1)}$};
\draw (4.5,1.5) node[draw=none,fill=none] {$l_{\gamma^{-1}(n)}$};
\draw (5,2.6) node[draw=none,fill=none] {$y_{\gamma^{-1}(n)}$};

    \draw[->] [ultra thick,blue] (2,2) to[out=45,in=135]  (5,2);
       \draw[->] [ultra thick,blue] (5,2) to[out=45,in=135]  (8,2);
       % \draw [ultra thick,blue] (3,2) to[out=-45,in=-135]  (7,2);
        \draw (8.7,2) node[draw=none,fill=none] {$y_{\gamma^{-2}(n)}$};
          \draw[<-] [ultra thick,blue] (2,2) to[out=-45,in=-135]  (5,0);
            \draw[<-] [ultra thick,blue] (5,2) to[out=-45,in=135]  (5,0);

        \draw[->] [ultra thick,blue] (5,0) to[out=-45,in=-135]  (8,2);
          \draw[->] [ultra thick,blue] (8,2) to[out=45,in=90]  (9.5,2);
            \draw[->] [ultra thick,blue] (0.5,2) to[out=90,in=135]  (2,2);
          
        \draw (5,-0.75) node[draw=none,fill=none] {$x$};
     %   \draw [decoration={
        %     markings,
           %  mark=at position 0.4 with \arrow{>}}
       %,postaction=decorate]
       %(A) to (B);
\end{tikzpicture}
\end{center}

\medskip
\subsection{General case | definition of $\mD$: ``Type I triples''}
\label{GCD}

We are now ready to define the class $\mD$.
Elements of $\mD$ are constructed through the following data: 
\begin{equation}
\mathcal{E} = (S,L,L',n,T)
\end{equation}
where 
\begin{itemize}
%\item
%$L$ and $L'$ are not conjugate (and in particular $L \ne L'$),
\item
$(S,L,L')$ is an EC-triple with FF (we call $(S,L,L')$ the {\em base} of $\mE$), 
\item
$n \in \mathbb{N}^\times$ is a positive integer, and 
\item
$T$ is any transitive subgroup of the symmetric group $\mathrm{S}_n$ acting naturally on $\{1,\ldots,n\}$. 
\end{itemize}

We construct the EC-triple $(S \wr T,L \wr T,{L'} \wr T)$ as above, which we also denote by $(\mE(S),\mE(L),\mE(L'))$. Note that 
$\vert S  \wr T\vert = \vert S\vert^n\vert T\vert$,   $\vert L \wr T\vert = \vert L\vert^n\vert T\vert$ and  $\vert {L'} \wr T\vert = \vert L'\vert^n\vert T\vert$.\\

If $L$ and $L'$ are not conjugate (in particular $L \ne L'$), then $L \wr T$ and $L' \wr T$ also are no conjugate. For, suppose the latter {\em are} conjugate: $(L \wr T)^g = L' \wr T$ with $g \in S \wr T$. Then as $(L \wr T) \cap S^n = L^n$ (where as usual we write $L^n$ for $L^n \times \{ \id\}$, etc.), $(L' \wr T) \cap S^n = {L'}^n$ and $S^n \unlhd S \wr T$, 
we have 
\begin{equation}
{(L^n)}^g = ((L \wr T) \cap S_n)^g = (L \wr T)^g \cap S_n^g = (L \wr T)^g \cap S_n = (L' \wr T) \cap S_n = {L'}^n.
\end{equation}
Now with $g = (\overline{u},\gamma)$, we obtain that
\begin{equation}
{(L^n)}^g = {(L^n \times \{\id\})}^g = (\gamma^{-1}(\overline{u}^{-1}L^n\overline{u})) \times \{\id\} = {(L')}^n \times \{\id\},
\end{equation}
and it easily follows that $L$ and $L'$ are also conjugate in $S$.\\

\medskip
\begin{remark}[Transitive groups]
{\rm
{\em Any} group acts transitively on some set: for instance, let it act on the left on itself, and one obtains a sharply transitive action.
}
\end{remark}

\medskip
\begin{remark}[Finiteness versus infinite]
{\rm
For physical means, one requires that the groups be finite, but everything works fine for infinite groups as well.
}
\end{remark}

\medskip
\subsection{Generic properties}

As we have seen, several of the ``naive'' product constructions generate examples which lose fundamental properties
of the initial base example, such as FF or MAX. We will show that our construction preserves them. This is most important for
the physical viewpoint.\\

\bt[Conservation of Properties]
Let $\mE = (S,L,L',n,T)$ be an element of $\mD$, and suppose $(S,L,L')$ is MAX.
Then $(\mE(S),\mE(L),\mE(L'))$ is an EC-triple which is FF, and also MAX.
\et

{\em Proof}.\quad
We first prove FF.
Suppose by way of contradiction that $N$ is a normal subgroup of $\mE(S)$ which is contained in $\mE(L)$. Then $N \cap S^{\times n}$ 
(where we write $S^{\times n}$ for $S^{\times n} \rtimes \{ \id\}$) is a normal subgroup of $S^{\times n}$ (as $S^{\times n}$ is a normal subgroup in $\mE(S)$). Since $N$ is contained in $\mE(L)$, 
it follows that $N \cap S^{\times n} = N \cap L^{\times n} \leq L^{\times n}$. For, an element of $N \cap S^{\times n}$ is of the form $(\overline{s},\id)$ with $\overline{s} \in S^{\times n}$, and as $N \leq L^{\times n} \rtimes T$, each element of $N$ has the form $(\overline{l},t)$ with $\overline{l} \in L^{\times n}$, $t \in T$. So the elements in $N \cap S^{\times n}$ are precisely the elements in $N \cap L^{\times n}$. 
By Principle \ref{sgr}, we conclude that $N \cap 
L^{\times n}$ must be trivial (if $(S,L,L')$ has FF, $(S^{\times n},L^{\times n},{L'}^{\times n})$ must have FF). Now let $(n,\gamma)$ be in $N \leq \mE(L)$; as $T \leq \mE(L)$, we have 
$(n,\id) \in \mE(L)$, which implies it to be an element of $L^{\times n}$. This is only possible if $n = \id$. So $N \leq T$.
Finally, let $(\id,\beta)$ be any element of $N$; then with $(r,\alpha)$ an arbitrary element of $\mE(S)$, we must have that 
\begin{equation}
(\id,\beta)^{(r,\alpha)} = (\alpha^{-1}(r^{-1}\beta(r)),\beta^\alpha)
\end{equation}
is an element of $N$. So $\alpha^{-1}(r^{-1}\beta(r)) = \id$, implying that $\beta(r) = r$ for all $r \in S^{\times n}$. 
It follows that $N$ is trivial, and that $(\mE(S),\mE(L),\mE(L'))$ has FF (since the same argument works for $\mE(L')$).\\

We now turn to MAX.
Suppose $\mE(L)$ is not a maximal subgroup of $\mE(S)$, and let $M$ be a proper subgroup of $\mE(S)$ that properly contains $\mE(L)$.  
Note that $T$ is contained in $M$ (up to a natural isomorphism $t \in T \mapsto (\id,t)$). Let $(m,\beta) \in M \setminus \mE(L)$. As $(\id,\beta^{-1}) \in M$, 
it follows that $(m,\id) \not\in L^{\times n}$, so with $m = (m_1,\ldots,m_n)$, some $m_i$ is not in $L_i \cong L$ (= the $i$th copy of $L$ 
in the Cartesian product $L^{\times n}$). So if $\pi_i: M \longrightarrow S_i \cong S$ is the projection of $M$ onto the $i$th copy
of $S$ in the product $S^{\times n}$, we have that $\pi_i(M)$ properly contains $L_i$, so that it must coincide with $S_i$ as $L_i$ is maximal
in $S_i$ by assumption. 
Now $L^{\times n} \rtimes \{ \id\} \leq M$, and if $\underline{M}$ is the projection of $M$ onto $S^{\times n} \rtimes \{\id\}$, then $L^{\times n} \rtimes \{\id\} \ne \underline{M}$. Consider the subgroup $L \times \{ \id\} \times \cdots \times \{ \id\}$ of $\underline{M}$. As $\pi_1(M) = S_1 = S$, we have that
\begin{equation}
\langle (L \times \{ \id\} \times \cdots \times \{ \id\})^{\underline{M}} \rangle = \langle L^S \rangle \times \{ \id\} \times \cdots \times \{\id\}. 
\end{equation}

But $L$ is maximal in $S$, so either $\langle L^S \rangle =  L$, and then $L$ is a normal subgroup of $S$, or $\langle L^S \rangle = S$. 
First suppose that we are in the latter case (for both the properties AC/EC), so that $S \times \{ \id\} \times \cdots \times \{ \id\} \leq \underline{M}$. 
As $T$ acts transitively on the components in $S^{\times n}$, and as $T$ is contained in $M$, 
it follows that $M$ coincides with $\mE(S)$. The same argument works for $\mE(L')$, so we have indeed proved that MAX 
is inherited.\\

Now let $L \unlhd S$. If $(S,L,L')$ is supposed to be AC, it follows as before that $L = L'$, contradiction. If $(S,L,L')$ is EC, as before it follows
that $L' \leq L$. On the other hand, as $L'$ is also supposed to be maximal in $S$, we must have  $L = L'$, contradiction. 
\eop \\

\subsection{Examples of unbounded index}

Since for all $m \in \mathbb{N}^{\times}$ we can find sharply transitive permutation groups $(T,X)$ with $\vert T\vert = \vert M\vert = m$ (e.g., $T$ is cyclic of order $m$ acting on itself by translation), we have the following theorem.

\bt[Unbounded index]
For each $N \in \mathbb{N}$, there exist EC-triples $(G,U,V)$ with MAX and FF, such that
\begin{equation}
\frac{\vert G\vert}{\vert U\vert} = \frac{\vert G\vert}{\vert V\vert} > N.
\end{equation}
\et
{\em Proof}.\quad
Follows immediately from the expressions of the sizes of the EC-triples $(S \wr T,L \wr T,{L'} \wr T)$ with $(S,L,L',n,T)$ in $\mD$ (by suitably adapting $n$ in function of $N$), and the remark preceding the theorem.
\eop \\

\subsection{Variation on the construction in \S\S \ref{GCD}: ``Type II triples''}

There is a second class of examples, $\widetilde{\mD}$, which we describe below. It is important to note that $S$ will be a {\em simple group} by assumption.

Elements of $\widetilde{\mD}$ are constructed through the following data: 
\begin{equation}
\mathcal{E} = (S,L,L',n,T)
\end{equation}
where 
\begin{itemize}
\item
$S$ is a simple group,
\item
$n \in \mathbb{N}^\times$ is a positive integer,  
\item
$L$ and $L'$ are diagonal subgroups of $S^n$,
\item
$(S^n,L,L')$ is an EC-triple (again called the {\em base} of $\mE$), and
\item
$T$ is any transitive subgroup of the symmetric group $\mathrm{S}_n$ acting naturally on $\{1,\ldots,n\}$, and both $L, L'$ are stabilized by $T$ in its natural action on $S^n$. 
\end{itemize}

We construct the triple $(S \wr T,L \rtimes T,{L'} \rtimes T)$ as above, which we also denote by $(\mE(S),\mE(L),\mE(L'))$. Note that in this construction, $\vert S  \wr T\vert = \vert S\vert^n\vert T\vert$,   $\vert L \rtimes T\vert = \vert L\vert\cdot \vert T\vert$ and  $\vert {L'} \rtimes T\vert = \vert L'\vert\cdot \vert T\vert$.\\

The fact that $(\mE(S),\mE(L),\mE(L'))$ is EC can be obtained as before.\\

Again, one can verify MAX and FF for these examples (note that for MAX, one does not have to assume MAX for $(S,L,L')$). We will quickly browse through the proofs, which are a bit different than for the original construction (due, a.o., to the fact that $S$ is simple).\\

\quad MAX.\quad One can derive MAX from O'Nan Scott theory (cf. the appendix of this paper), but we will describe the idea for the 
sake of the reader.
Without loss of generality, suppose
\begin{equation}
L := \{ (s,s,\ldots,s) \vert s \in S\},
\end{equation}
and denote $S^n$ by $N$. Suppose by way of contradiction that $M$ properly contains $L \rtimes T$, $M \ne N \rtimes T$. Then $L \ne M \cap N \ne N$. Let $R \geq M \cap N$ be maximal in $N$. Before proceeding, recall that a group $C$ is maximal in a direct product $D^r$ (of $r$ copies of $D$, $r \geq 2$), if either it is a normal subgroup of prime index, or $C = \pi_{i,j}^{-1}(E)$, where $E$ is maximal in $D_i \times D_j$, and $\pi_{i,j}$ is the projection of $D^r$ on the $i$-th and $j$-th component in $D^r$.
We call a maximal subgroup of $D^r$ {\em standard} if it is a  direct product of a maximal subgroup in $D$ with $r - 1$ copies of $D$. If 
$R$ is a normal subgroup of prime index in $N$, then $R$ does not contain some $S_i$ (seen as a subgroup of $N$), so $R \cong S^{m - 1} \times (S_i \cap R)$. By the transitive action of $M$ on the components, it follows that $R = N$. If $R$ is not normal of prime index, $R = \pi_{i,j}^{-1}(M')$ for some maximal subgroup $M'$  in $S_i \times S_j$. If $m \geq 3$, then $R \cong S^{m - 2} \times M'$, and the transitivity argument of above leads to $R = N = M\cap N$, and hence $M = S^n \rtimes T$. If $m = 2$, then we use the fact that a group $F$ is simple if and only if the diagonal group is maximal in $F \times F$, to conclude that $L$ is already maximal in $S \times S = N$, so that $M \cap N = N$, and hence $M = S^n \rtimes T$.\\

\medskip
\quad FF.\quad
Suppose $K$ is a normal subgroup of $N \rtimes T$ which is a subgroup of $L \rtimes T$. First suppose that $(\overline{k},\id) \in K$ (with 
$\overline{k} = (e,e,\ldots,e) \in L^n$). Then for any $(\overline{v},\id)$ in $N \rtimes T$, $\overline{v} = (v_1,\ldots,v_n)$, we have that $(\overline{k},\id)^{(\overline{v},\id)} = (e^{v_1},\ldots,e^{v_n},\id) \in K \leq L \rtimes \{\id\}$. This is obviously not possible: take, for instance, $v_1 = \id$ and 
$v_2 \not\in C_S(e)$ (which is possible as $S$ is simple). (For $n = 1$, there is nothing to prove by assumption.)\\

It follows that $K \cap (L \rtimes \{\id\}) = \{\id\}$.  Now let $(\overline{k},\gamma) \in K$ (with 
$\overline{k} = (e,e,\ldots,e) \in L$ and $\gamma \ne \id$). Then for any $(\overline{v},\id)$ in $N \rtimes T$, $\overline{v} = (v_1,\ldots,v_n)$, we have that $\varphi := (\overline{k},\gamma)^{(\overline{v},\id)} = (\overline{v}^{-1}\overline{k}\gamma(\overline{v}),\gamma)$. Let $\overline{v} := (s,\id,\ldots,\id)$.
If $n \geq 3$, then at least on entry of $\varphi$ is simply $e$, while the others have the form $s^{-1}es, s^{-1}e, es$ or $e$, and not only the latter form occurs. If an entry has the form $s^{-1}es$, let again $s \not\in C_S(e)$; if the entry has the form $s^{-1}e$ or $es$, then let $s \ne \id$. It follows that we can find $\overline{v}$ such that $(\overline{k},\gamma)^{(\overline{v},\id)}$ is not in $K$, that is, such that $K$ is not a normal subgroup of $N \rtimes T$ inside $L \rtimes T$.\\

We have shown that FF holds (as the same works for $L'$).

\medskip
\bt[Conservation of Properties]
Let $\mE = (S,L,L',n,T)$ be an element of $\widetilde{\mD}$.
Then $(\mE(S),\mE(L),\mE(L'))$ is an EC-triple which is FF, and also MAX.\eop
\et

\medskip
\subsection{Second variation on the construction in \S\S \ref{GCD}: ``Type III triples''}

There is a second class of examples, denoted $\widehat{\mD}$, which we describe below. It is important to note that $S$ again will be a {\em simple group} by assumption. These examples will be a mixture of the two first constructions.

Elements of $\widehat{\mD}$ are constructed through the following data: 
\begin{equation}
\mathcal{E} = (S,L,L',(l,k),T)
\end{equation}
where 
\begin{itemize}
\item
$S$ is a simple group,
\item
$l, k \in \mathbb{N}^\times$ are positive integers both different from $1$,
\item
$L$ and $L'$ are diagonal copies of $S$ in (groups isomorphic to) $S^k$,
\item
$(S^k,L,L')$ is an EC-triple with FF (again called the {\em base} of $\mE$),  and
\item
$T$ is any transitive subgroup of $\mathrm{S}_{lk}$ acting naturally on $\{1,\ldots,lk\}$, by a block system with block $\{1,\ldots,k\}$, and
both $L, L'$ are stabilized by $T_{\{1,\ldots,k\}}$ in its natural action on $S^k$.
\end{itemize}

We construct a triple $(S \wr T,L^l \rtimes {T},{L'}^l \rtimes {T})$ as above, which we also denote by $(\mE(S),\mE(L),\mE(L'))$. Note that 
$\vert S  \wr T\vert = \vert S\vert^{lk}\vert T\vert$,   $\vert L^l \rtimes T\vert = \vert S\vert^l\vert T\vert$ and  $\vert {L'}^l \rtimes T\vert = \vert S\vert^l\vert T\vert$.\\

The fact that $(\mE(S),\mE(L),\mE(L'))$ is EC can again be obtained as before.\\

The proof of the following result is left to the reader.

\medskip
\bt[Conservation of Properties]
Let $\mE = (S,L,L',(l,k),T)$ be an element of $\widehat{\mD}$.
Then $(\mE(S),\mE(L),\mE(L'))$ is an EC-triple which is FF, and also MAX.\eop \\
\et

\begin{remark}
\label{remblock}
{\rm
Note that we suppose that $L$ and $L'$ are diagonal copies of $S$ in groups isomorphic to $S^k$. After applying the appropriate isomorphisms, we can then indeed without loss of generality formulate the next property for $\widehat{\mD}$ as ``$(S^k,L,L')$ is an EC-triple with FF'' (not having changed the notation for $L$ and $L'$ for the sake of convenience). In exactly the same vein, after considering the last defining property for  
$\widehat{\mD}$ (namely, the one concerning the action of $T$), $L$ and $L'$ could be defined as diagonal copies of $S$ inside respectively 
$S_{u_1} \times \cdots \times S_{u_k}$ and $S_{u_1'} \times \cdots \times S_{u_k'}$, where $\{u_1,\ldots,u_k \}$ and $\{u_1',\ldots,u_k' \}$ are members of the same block system described above.
}
\end{remark}

 \section{Classifying planar isospectral drums}
 
 Let $(G,U,V)$ be a (finite) EC-triple with MAX and FF. 
 
 We will determine the structure of such an example, imposing more properties on the way which
 are satisfied for {\em any pair of transplantable planar isospectral drums}. Eventually, we will observe that 
finite (almost) simple groups  yield the right triples.

To some point, we will follow some standard steps in the proof of the O'Nan-Scott Theorem (cf. the appendix). But especially the physical constraints will eventually shape the classification (in \S\S \ref{case3} and \S\S \ref{case4}).\\

 We will use the following fact throughout (see the appendix).\\
 
 \textsc{Fact}.\quad 
 {\em If $X$ is a maximal subgroup of $Y$ and the action of $Y$ by left multiplication on the left cosets of $X$
 is faithful, then for a hypothetical nontrivial normal subgroup $Z$ of $Y$, we have that $Y = XZ$. In other words, $Z$ acts transitively on $\{ lX\ \vert\ l \in Y \}$.}\\

 \subsection{Case I: $G$ is simple}
 
 If $G$ is simple, we are in a special case of our construction, Type I, namely the case with data $\mE = (G = S,U = L,V = L',1,\{\id\})$.
 
 Note that all the known examples of isospectral pairs fall under Case 1 (in each of these examples, $G \cong \PSL_\ell(q)$ for small values
 of $\ell$ and $q$, and these groups are simple).

 \subsection{Case II: $G$ is not simple | first case}
 
 As $U$ is maximal, $G$ acts primitively on the left cosets of $U$ by left multiplication (and the same can be said about $V$).
 The fact that $G$ is not simple implies that there are nontrivial minimal normal subgroups. But, there are at most two (cf. the appendix), and 
 detailed information is known about them. First we suppose there is only {\em one} minimal normal subgroup $N$, and that it is not abelian.
 The fact that there is one and only one nontrivial minimal normal subgroup implies that the action of $G$ on the left coset space $G/U$ is not sharply transitive (cf. the appendix), that is, $U \cap N \ne \{ \id\}$.

 By Appendix \ref{Wil}, $N$ is (isomorphic to) a direct product of ($m$) copies of the same nonabelian simple group, say $S$. 
 %Also, $C_G(N) = \{\id\}$.
  Now $G$ acts by conjugation as an automorphism group on $N$. 
 %, and as $C_G(N) = \{\id\}$, the action is faithful. \\
 
 \quad{\bf First put $m = 1$, so that $N = S$ is simple.}\quad 
 Then $G$ is an {\em almost 
 simple group} (by definition), since $S \leq G \leq \Aut(S)$. So this falls under the data $\mE = (G = S,U = L,V = L',1,\{\id\})$ for the Type I construction.\\

 \quad{\bf Now let $m \geq 2$.}\quad 
 Write $N$ as
 \begin{equation}
 N = S_1 \times \cdots \times S_m,
 \end{equation}
 where each $S_i \cong S$ is simple, and let $\mS = \{ S_1,\ldots,S_m\}$.
 
 First note that as $G = UN = VN$, it follows that 
 \begin{equation}
 G/N = UN/N \cong U/(U \cap N) \cong V/(V \cap N),
 \end{equation}
 that is, $U$ and $V$ act in the same way on $\mS$ as $G$ does (they ``contain'' the complete $G$-action). By Appendix \ref{Wil}, the 
 action is transitive.

 By O'Nan-Scott theory (see for instance the excellent lecture notes \cite{RAW} for more details besides those in the appendix), there are two cases (with some subdivision) to consider (for both $U$ and $V$ | we only state them for $U$):
 \begin{itemize}
 \item[A]
{\em  The natural maps $\pi_i: U \mapsto S_i$ are not surjective.}\quad  This implies that $U \cap N  = M_1 \times \cdots \times M_m$, where each $M_i \cong M$ is a maximal subgroup in $S$. It follows that $\vert U \cap N \vert = \vert M\vert^m$.
 \item[B]
  {\em The natural maps $\pi_i: U \mapsto S_i$ {\em are} surjective.}\quad  There are two subcases.
  \begin{itemize}
  \item[B.1]
  $U \cap N$ is a diagonal copy $\mathrm{diag}(S)$ of $S$ in $N$. In particular, $\vert U \cap N \vert = \vert S\vert$.
  \item[B.2]
  $U$ acts on $\{1,\ldots,m\}$ through blocks $\{u_1,\ldots,u_{k|}\}^U$, where $\chi(U) := \{ u_1,\ldots,u_k\} \subset \{ 1,\ldots,m\}$, and 
  $k$ divides $m$, and
  $U \cap (S_{u_1} \times \cdots \times S_{u_k})$ is a diagonal copy $\mathrm{diag}(S)$ of $S$. (So in this case, $\chi(U)$ is a block in a system of blocks for the action of $G$ on $\{ 1,\ldots,m\}$.) In particular, 
  $\vert U \cap N \vert = \vert S \vert^{m/k}$.
  \end{itemize}
 \end{itemize}
 
 If we are in (A), then  $U \cong M \wr T$, where $T$ is isomorphic to the subgroup of $\mathrm{S}_m$ which is induced on $\mS$ by $G$, and also $U$. Also, $G \cong S \wr T = N \rtimes T$. \\ 

If we are in (B.1), $U \cong \mathrm{diag}(S) \rtimes T$, with $T$ as above. Also, $G \cong S \wr T = N \rtimes T$.\\

If we are in (B.2), $U \cong \mathrm{diag}(S)^{m/k} \rtimes {T}$, and also, $G \cong S \wr T = N \rtimes T$. It should be noted that   $\mathrm{diag}(S) \leq N$. \\ 

%we have an action of any $S_i$ on the left coset space $S_i/M_i$ contained in $\mathbf{S}_r$, with $r = [S : M]$
 
% Let $K$ be the subgroup of $G$ which fixes each element of $\mS := \{ S_1,\ldots,S_m\}$ when acting by conjugation (so $K$ is by definition
 %the {\em kernel} of this action). Then $K$ is a normal subgroup of $G$ which contains $N$, and $(G/K,\mS)$ is a faithful permutation group.
 %  Put $G/K = B$, with $B$ isomorphic to  a transitive subgroup of $\mathrm{Sym}(\mS) = S_m$. As obviously $K \leq 
 %\mathrm{Aut}(\mS) \wr S_m$, $K \cong A_1 \times \cdots \times A_m$, with each $A_i \cong A$, and $S \leq A \leq \mathrm{Aut}(S)$. 
 %So 
 %\begin{equation}
 %G \cong K \rtimes B \cong A \wr B.
 %\end{equation} 
 %So $(G,U,V)$ has the form $(A \wr B,(K \cap U) \rtimes B,(K \cap V) \rtimes B)$.\\

 Before finishing this part of the classification, we need the following result. First note that in (A), (B.1) and (B.2), $(G,U,V)$ already has the form
 as in the Type I, II, III triples (respectively) defined before. In the same way as for these triples, we introduce the ``base triple.'' 
 
 \begin{theorem}
 Let $(G,U,V)$ be EC. Let $(A,B,C)$ be the base triple in each of (A), (B.1), (B.2).
 \begin{itemize}
 \item[{\rm (1)}]
  We have that $(A,B,C)$ is EC. 
  \item[{\rm (2)}]
  If $U$ and $V$ are not conjugate, then $B$ and $C$ are not conjugate.
  \item[{\rm (3)}]
  If $\vert U\vert = \vert V\vert$, then $\vert A\vert = \vert B\vert$, and hence $U$ and $V$ {\em both} are constructed through  the {\em same} process (A), (B.1) or (B.2).\\
  In particular, if PAIR is satisfied, the same conclusion holds.
  \end{itemize}
 \end{theorem}
 
 {\em Proof}.\quad
 (1)\quad
 Take any $(\overline{u},t) \in U$; then by EC there is a $(\overline{n},r) \in G$ and $(\overline{v},t') \in V$ such that
 \begin{equation}
 (\overline{u},t)^{(\overline{n},r)} = (\overline{v},t').
 \end{equation}
 In particular, the formula holds when $t = \id$, $t' = \id$, and we then have
 \begin{equation}
  (\overline{u},\id)^{(\overline{n},r)} = (r^{-1}(\overline{n}^{-1}\overline{u} \overline{n}),\id) = (\overline{v},\id).
 \end{equation}
 So for all $\overline{u} = ((b_1,\ldots,b_m),\id)$ in $U$ (each $b_i$ taken in $B$), there is some $(\overline{n},\id) = ((a_1,\ldots,a_m),r)$ in $G$ (each $a_i$ taken in $A$)
 and $(\overline{v},\id) = ((c_1,\ldots,c_m),\id)$ in $V$ (each $c_i$ taken in $C$), such that 
 \begin{equation}
 ((b_1^{a_1},\ldots,b_m^{a_m}),\id) = (r(c_1,\ldots,c_m),\id).
 \end{equation}
 The element $(r(c_1,\ldots,c_m),\id)$ is an element of $V$, as $V$ also contains the $G$-action.
  The EC-property easily follows for $(A,B,C)$ in each of the cases (A), (B.1) and (B.2).\\
  
  (2)\quad 
  Was obtained before.\\
  
  (3)\quad
  Follows immediately.
 \eop \\
 
 \begin{remark}{\rm
 Consider the triple $(G,U,V)$ in B.2 (so both $(G,U)$ and $(G,V)$ are assumed to fall under B.2). Then $(G,U) = (S\wr T,U\cong S^l \rtimes T)$ and $(G,V) = (S \wr T,V\cong S^l \rtimes T)$.
 Also, $G$ acts on $N$, and on $\mS = \{S_1,\ldots,S_m\}$ by conjugation, and $T$ also does. Besides that, $T$ preserves the 
 block systems of $U$ and $V$. If $u \in U$, there is some $g \in G$ and $v \in V$ such that $u^g = v$ by EC, so it follows that $U$ and $V$ live in the same block system. So we are indeed dealing with Type III triples (that is, not only $(G,U)$ and $(G,V)$ have the desired form, but $(G,U,V)$ as well, as it should have).
 }
 \end{remark}

 \medskip
 From this point on, the {\em physical properties} (FF and INV) will be the deciding factors in the classification.
 
 \medskip
 \subsection{Case III: $G$ is not simple | second case}
 \label{case3}
 
 Now suppose $G$ has precisely two minimal normal subgroups $N$ and $N'$; then by Appendix \ref{Wil}, $N = C_H(N')$, $N' = C_H(N)$,
 $N \cap N'$ and $N \cap N' = \{\id\}$, so that $M := N \times N'$ is a normal (characteristic) subgroup of $G$.
 %As $N \cap N' = \{\id\}$, $G$ acts by conjugation as a faithful automorphism group on $M$. 
 Both $N$ and $N'$ act sharply transitively 
 on the left cosets of $U$ and $V$, so that $U \cap N = U \cap N' = V \cap N = V\cap N' = \{\id\}$, and
 \begin{equation}
 G =  UN = UN' = VN = VN'.
 \end{equation}
 
 In particular, $\vert U\vert = \vert V\vert$. \\ 
 
 By INV, there is an involution $\gamma$ in $G$ which has more than $\vert \mX\vert/3$ fixed points in $\mX$. If $x$ and $y$ are 
 any two of these fixed points, and
 $n \in N$ is such that $x^n = y$, then as $N$ acts sharply transitively on $\mX$ and as it is a normal subgroup of $G$, it follows that $n$ and $\gamma$ commute (as the commutator $[n,\gamma]$ fixes $y$ and is contained in $N$). So 
 \begin{equation}
 \mathrm{Fix}(\gamma) = \vert C_N(\gamma)\vert > \frac{\vert \mX\vert}{3} = \frac{\vert N\vert}{3}.
 \end{equation}
 
 So only the possibilities $\vert N \vert/2$ and $\vert N\vert$ come out for $\vert C_N(\gamma)\vert$. If $N = C_N(\gamma)$, it follows that 
 $\gamma$ fixes all points of $\mX$, a property which violates FF (and also INV). 
 
 Now turn to the case $\vert C_N(\gamma)\vert = \vert N\vert/2$.
 By Appendix \ref{Wil}, $N$ is  (isomorphic to) a direct product of ($m$) copies of the same nonabelian simple group, $S$. 
Write 
 \begin{equation}
 N \cong S_1 \times \cdots \times S_{m},
 \end{equation}
 each $S_j$ being isomorphic to $S$.

Some $S$-factor $\widetilde{S}$ of $N$ is not contained in $C_N(\gamma)$ (as otherwise $N = C_N(\gamma)$), and then $\widetilde{S} \cap C_N(\gamma)$ is a subgroup of index $2$ in $\widetilde{S}$. So if $\vert S\vert > 2$, $\widetilde{S} \cap C_{N}(\gamma)$ is a nontrivial
normal subgroup of $\widetilde{S}$, contradicting its simplicity. So $\vert S \vert = 2$, meaning that $N$ is elementary abelian, contradiction. 
(This is in fact Case IV | see \S\S \ref{case4}.)\\

 %Put $\mS = \{ S_1,\ldots,S_{m} \}$, and let, as above, $K$ be the kernel of the action of $G$ on $\mS$. Remark that $G$ acts transitively on $\mS$ %(cf. Appendix \ref{Wil}).
 %Then $K$ is a normal subgroup containing $N$, so that $G = UK = VK$, 
 %and as above we have
 %\begin{equation}
 %G/K \cong U/(U \cap K) \cong V/(V \cap K). 
 %\end{equation}
 %So again we have that $U$ and $V$ ``contain'' the action of $G$ on $\mS$.
 %Now repeat the argument of the second part of Case II.\\

 \medskip
 \subsection{Case IV: $G$ is not simple | abelian case}
 \label{case4}
 
 Suppose $N$ is abelian. In that case, $N$ acts sharply transitively on the left coset spaces $G/U$ and $G/V$, and 
 $N$ is elementary abelian (cf. the appendix to this paper):
 \begin{equation}
 N \cong C_p \times C_p \times \cdots \times C_p,
 \end{equation}
 where $C_p$ is the cyclic group of prime order $p$, and where we have taken $m$ copies on the right-hand side.
So we can identify $N$ as the translation group  of the affine space $\AG(m,p) =: \bA$, which is itself identified with, say, $G/U$. 
(So the elements of $G/U$ correspond to the points of $\AG(m,p)$.) After this identification, $G$ acts as a faithful automorphism 
group of $\bA$, and it is isomorphic to a subgroup of $\AGL_{m + 1}(p)$. 

We will now invoke INV to solve the last piece of the puzzle. As INV assumes an identity involving three involutions which are contained in $G$,
we need extra knowledge of how involutions act on affine spaces. So we will discuss the different types of involutions that can occur in the automorphism group of a finite
{\em projective} space $\PG(n,q)$, cf.  \cite{BS}. The reader can deduct the different types of involutions for {\em affine} spaces from this result.

\begin{itemize}
\item[{\bf Baer involutions}]
A {\em Baer involution} is an involution which is not contained in the linear automorphism group of the space, so that $q$ is a square, and it fixes an 
$n$-dimensional subspace over $\mathbb{F}_{\sqrt{q}}$ pointwise. 
\item[{\bf Even characteristic}]
If $q$ is even, and $\theta$ is an involution which is not of Baer type, $\theta$ must fix an $m$-dimensional subspace of 
$\PG(n,q)$ pointwise, with $1 \leq m \leq n \leq 2m + 1$.
In fact, to avoid trivialities, one assumes that $m \leq n - 1$.
\item[{\bf Odd characteristic}]
If $\theta$ is a linear involution of $\PG(n,q)$, $q$ odd, the set of fixed points is the union of
two disjoint complementary subspaces. Denote these by $\PG(k,q)$ and $\PG(n - k -1,q)$, $k \geq n - k - 1 > -1$.
\item[{\bf Other involutions}]
Those without fixed points.
\end{itemize}

The argument below can now to some extent be taken from the combined work of Giraud \cite{Giraud} and the author \cite{KTI,KTII,KTIII}.

By INV we want to consider triples $(\mathcal{A},\{\theta^{(i)}\},r)$, 
where $\mathcal{A}$ is our affine space of dimension $m \geq 2$ over $\F_p$,
and $\{\theta^{(i)}\}$ a set of $r$ nontrivial involutory automorphisms of $\mathcal{A}$, satisfying
\begin{equation}
\label{eq12}
 r(\vert \mathcal{A}\vert) - \sum_{j=1}^r\mbox{Fix}(\theta^{(j)}) = 2(\vert \mathcal{A}\vert - 1),
\end{equation}
for some natural number $r \geq 3$. So we consider
\begin{equation}
(r - 2)p^m + 2 = \sum_{j=1}^r\mbox{Fix}(\theta^{(j)}).         
\end{equation}

As $p$ is a prime, Baer involutions do not occur.
Since an automorphic involution fixes at most $p^{m - 1} + 1$ points of $\mathcal{A}$, we have $r = 3$. It  
is also clear that $p \leq 3$. 

For $p = 3$, the only solution is $\{\mathrm{Fix}(\theta^{(1)}),\mathrm{Fix}(\theta^{(2)}),\mathrm{Fix}(\theta^{(3)})\}
= \{p^{m - 1} + 1,p^{m - 1} + 1,p^{m - 1}\} =  \{3^{m - 1} + 1,3^{m - 1} + 1,3^{m - 1}\}$. Without loss of generality, we suppose $\theta^{(1)} =: \theta$ has $3^{m - 1} + 1$ fixed points in $\mA$, so that $\theta$ fixes an affine subspace $\mA_{\theta}$ of dimension $m - 1$ pointwise, and one additional affine point $z$.  Let $\Delta$ be the hyperplane at infinity of $\mA$. Below, if $\alpha \subseteq \mA$ is an affine subspace of $\mA$, by $\overline{\alpha}$ we denote the projective completion of $\alpha$, which is a projective subspace of $\overline{\mA}$. Also, if $\beta$ is an automorphism of $\mA$ (i.e., an element of $\AGL_{m + 1}(p)$), then $\overline{\beta}$ is the corresponding automorphism of $\overline{\mA}$.
First suppose $m \geq 2$, and let $y$ be any point of $\mA$ contained in $\Delta$, but not contained in the projective completion of $\mA_{\theta}$. Then in the projective completion of $\mA$, we have that $zy \cap \overline{\mA_{\theta}} = \{y'\}$ is a point which is not in $\Delta$. As $\overline{\theta}$ fixes $z$ and $y'$ but not $y$, it follows that $\Delta$ cannot be fixed by $\overline{\theta}$, contradiction by the definition of $\theta$. It follows that $m = 1$. In that case, each $\theta^{(i)}$ fixes precisely one point of $\mA$, and whence this case cannot occur.\\

Now let $p = 2$. Then
\begin{equation}
2^m + 2 = \sum_{j=1}^r\mbox{Fix}(\theta^{(j)}).         
\end{equation}

As the left-hand side is divisible by $2$ but not by $4$, some $\theta^{(i)}$ must fix precisely two points of $\mA$ (we are working in even characteristic!). So $(m - 1)/2 \leq 1$, and hence $m \leq 3$. If $m = 1$, it is easily seen that there are no solutions of the INV-equation. If $m \in \{2,3\}$, the only numerical solution is given by $\{2,2^{m - 1},2^{m - 1}\}$.\\

All AC-triples with these numerical properties are known | see e.g. \cite{OS}. We have not performed a more detailed analysis to find EC-triples with these data.\\

\medskip
\subsection{On AC-triples and planar examples}
\label{planar} 
 
One expects triples $(G,U,V)$ that yield elusive (counter) $\mathbb{R}^2$-examples to fall inside classes I, II | first case, III or IV; in class II | second case (so with $m \geq 2$), a similar analysis as in the affine case after invoking INV  might lead to a very restrictive list (read: empty) of arithmetic possibilities. If that were true, one would have that an AC-triple $(G,U,V)$ with MAX, FF, PAIR and INV (that is, an irreducible Gassman-Sunada triple which yields planar isospectral non-isometric billliards) is of {\em only two types}:

\begin{itemize}
\item[{\bf PLAN-1}]
of ``almost simple type,''
\item[{\bf PLAN-2}]
or known.
\end{itemize}    

The author hopes to perform this analysis soon.

\newpage
\section{Conclusion}
We started this paper with a detailed explanation of the connection between the notions of ``transplantability'' and ``almost conjugacy,'' which are crucial in the construction theory of planar isospectral billiards. We then discussed the related notions {\em AC-triple} (or {\em Gassman-Sunada triple}) and {\em EC-triple} (which naturally generalizes AC-triples), in the context of isospectral drumheads (in any dimension).  Next, we explored some ``naive'' constructions of EC-triples, starting from one given EC-triple, through direct products and ``adding kernels.'' Then we 
introduced four properties | being FF, MAX, PAIR, INV | inspired by physical properties which ``irreducible'' drumheads should have on the one hand, and which planar examples have, on the other.

We introduced an extremely general set of construction procedures of EC-triples  with FF and MAX (starting from a given one), which due to MAX, is connected to the O'Nan-Scott Theorem of finite Group Theory, which involves finite simple groups. The essential reason why the procedures work can be situated in properties of underlying combinatorial point-line geometries which translate EC.

In a final stage, we classified EC-triples, by starting from EC-triples with FF, MAX, PAIR and INV, and showed that, indeed, any such example is constructed through one of our procedures. In particular, all the known planar counter examples to Kac's question ``Can one hear the shape of a drum'' arise.

In such a way, we have obtained a deep connection between irreducible (possibly higher-dimensional) drumheads and finite simple groups that offers a new direction in the theory.

\newpage
\bigskip
\section{Appendix: The O'Nan-Scott theorem}
\label{Wil}

Let $H \leq \Sym(\Omega)$ be a finite primitive permutation group, let $M$ be a minimal normal subgroup of $H$,
and let $z \in \Omega$. 
%Then $M$ is a characteristically simple group, isomorphic to a direct product of 
%pairwise non-isomorphic finite simple groups.

Below, $\soc(Y)$, where $Y$ is a group, denotes the {\em socle} of $Y$ | the normal subgroup generated by the minimal normal subgroups of $Y$. Also,
$\hol(Y)$ denotes the {\em holomorph} of $Y$ | a group which is isomorphic to $Y \rtimes \Aut(Y)$.

The O'Nan-Scott theorem distinguishes the following cases:\\

\begin{itemize}
\item[(HA)]
$M$ is elementary abelian and regular, $C_H(M)= M$, and $H \leq \hol(M)$; in this case $M$ is the translation group of some finite affine space 
$\mathbf{AG}(l,p)$ ($p$ a prime and $l$ some positive integer), and $H \leq \mathbf{GL}_l(p)$.\\
\item[(HS)]
$M$ is non-abelian, simple and regular, $\soc(H) = M \times C_H(M) \cong M\times M$, and $H \leq \hol(M)$;\\
\item[(HC)]
$M$ is non-abelian, non-simple and regular, $\soc(H) = M \times C_H(M) \cong M\times M$, and $H \leq \hol(M)$;\\
\item[(SD)]
$M$ is non-abelian and non-simple, $M_z$ is a simple subdirect product of $M$ and $C_H(M) = \{\id\}$;\\
\item[(CD)]
$M$ is non-abelian and non-simple, $M_z$ is a non-simple subdirect product of $M$ and $C_H(M) = \{\id\}$;\\
\item[(PA)]
$M$ is non-abelian and non-simple, $M_z$ is not a subdirect subgroup of $M$ and $M_z \ne \{\id\}$; $C_H(M) = \{\id\}$;\\
\item[(AS)]
$M$ is non-abelian and simple, $C_H(M) = \{\id\}$ and $H$ is an almost simple group;\\
\item[(TW)]
$M$ is non-abelian and non-simple, $M_z = \{\id\}$ and $C_H(M) = \{\id\}$. Also, there is a (non-abelian) simple group $T$ such that
$M = T \times \cdots \times T$ ($m$ factors) (and so $\vert M \vert = \vert T\vert^m$).\\
\end{itemize}

\bigskip
The theorem has a simple proof.
Below, we list some fundamental lemmas which are used to prove O'Nan-Scott, and which were frequently used in the body of this paper. All of them can be found in \cite{RAW}.

\bo
Every normal subgroup of a primitive group is transitive.
\eo

\bo
If $K$ is a characteristically simple group, then it is a direct product of isomorphic simple groups. In particular, this applies to minimal normal subgroups of any finite group.
\eo

\bo
If $H$ is a primitive group, and $N$ is a minimal normal subgroup of $H$, then $H = KN$ for any point stabilizer $K$.
\eo

Based on the previous observations, we consider the following, more transparent form of O'Nan-Scott, which summarizes the previous version in four main classes.\\

Let $H$ be primitive (in its action on the set $\mX$), and let $N$ be a minimal normal subgroup of $H$. Let $x \in \mX$, and let $H_x$ be the stabilizer of $x$ in $H$.

\begin{itemize}
\item[{\bf I}]
If any such $N$ is trivial, $H$ is simple.
\end{itemize}

If no such $N$ is abelian, there are two cases to consider.
\begin{itemize}
\item[{\bf II}]
{\em $N$ is unique.} If $N$ is simple, $H$ is almost simple. If $N$ is not simple, $N$ is a direct product of ($r$) mutually isomorphic simple groups $\cong P$, on which $H$ acts transitively by conjugation. One distinguishes three cases: (a) $H_x \cap N = W^r \leq P^r$ is a direct product of ($r$) mutually isomorphic groups ($W \leq P$); (b) $H_x \cap N$ is a diagonal copy of $P$; (c) there is a block system for $H$ in its action on $\{ P_1,\ldots,P_r\}$, where $N = P_1 \times \cdots \times P_r$, and a block $B$ such that  $(H_x \cap N) \cap B$ is a diagonal copy of $P$. (Here, we see $B$ as the subgroup of $N$ naturally defined by its elements.)
\item[{\bf III}]
{\em $N$ is not unique.} In that case, there is a second minimal normal subgroup $N'$, and $N \cap N' = \{\id\}$, $C_H(N) = N'$, $C_H(N') = N$, so that 
$N \times N'$ is a normal subgroup of $H$. Also, $N$ and $N'$ have the same properties as $N$ in (II), and there is an outer automorphism of $H$ which acts as an involution on $\{ N,N'\}$. Finally, both $N$ and $N'$ have the same sharply transitive action on $\mX$.
\end{itemize}

Finally, there is the abelian case.
\begin{itemize}
\item[{\bf IV}]
If $N$ is abelian, it is elementary abelian and sharply transitive on $\mX$, and $H$ is an affine group.
\end{itemize}

For more information, we refer the reader to \cite{RAW}.

\end{document}